\input amstex

%

\def\next{AMS-SEKR}\ifx\styname\next \endinput\fi
\catcode`\@=11
\def\styname{AMS-SEKR}
\def\styversion{2.0}
{\W@{}\W@{\styname.STY - Version \styversion}\W@{}}
\hyphenation{acad-e-my acad-e-mies af-ter-thought anom-aly anom-alies
an-ti-deriv-a-tive an-tin-o-my an-tin-o-mies apoth-e-o-ses apoth-e-o-sis
ap-pen-dix ar-che-typ-al as-sign-a-ble as-sist-ant-ship as-ymp-tot-ic
asyn-chro-nous at-trib-uted at-trib-ut-able bank-rupt bank-rupt-cy
bi-dif-fer-en-tial blue-print busier busiest cat-a-stroph-ic
cat-a-stroph-i-cally con-gress cross-hatched data-base de-fin-i-tive
de-riv-a-tive dis-trib-ute dri-ver dri-vers eco-nom-ics econ-o-mist
elit-ist equi-vari-ant ex-quis-ite ex-tra-or-di-nary flow-chart
for-mi-da-ble forth-right friv-o-lous ge-o-des-ic ge-o-det-ic geo-met-ric
griev-ance griev-ous griev-ous-ly hexa-dec-i-mal ho-lo-no-my ho-mo-thetic
ideals idio-syn-crasy in-fin-ite-ly in-fin-i-tes-i-mal ir-rev-o-ca-ble
key-stroke lam-en-ta-ble light-weight mal-a-prop-ism man-u-script
mar-gin-al meta-bol-ic me-tab-o-lism meta-lan-guage me-trop-o-lis
met-ro-pol-i-tan mi-nut-est mol-e-cule mono-chrome mono-pole mo-nop-oly
mono-spline mo-not-o-nous mul-ti-fac-eted mul-ti-plic-able non-euclid-ean
non-iso-mor-phic non-smooth par-a-digm par-a-bol-ic pa-rab-o-loid
pa-ram-e-trize para-mount pen-ta-gon phe-nom-e-non post-script pre-am-ble
pro-ce-dur-al pro-hib-i-tive pro-hib-i-tive-ly pseu-do-dif-fer-en-tial
pseu-do-fi-nite pseu-do-nym qua-drat-ics quad-ra-ture qua-si-smooth
qua-si-sta-tion-ary qua-si-tri-an-gu-lar quin-tes-sence quin-tes-sen-tial
re-arrange-ment rec-tan-gle ret-ri-bu-tion retro-fit retro-fit-ted
right-eous right-eous-ness ro-bot ro-bot-ics sched-ul-ing se-mes-ter
semi-def-i-nite semi-ho-mo-thet-ic set-up se-vere-ly side-step sov-er-eign
spe-cious spher-oid spher-oid-al star-tling star-tling-ly
sta-tis-tics sto-chas-tic straight-est strange-ness strat-a-gem strong-hold
sum-ma-ble symp-to-matic syn-chro-nous topo-graph-i-cal tra-vers-a-ble
tra-ver-sal tra-ver-sals treach-ery turn-around un-at-tached un-err-ing-ly
white-space wide-spread wing-spread wretch-ed wretch-ed-ly Brown-ian
Eng-lish Euler-ian Feb-ru-ary Gauss-ian Grothen-dieck Hamil-ton-ian
Her-mit-ian Jan-u-ary Japan-ese Kor-te-weg Le-gendre Lip-schitz
Lip-schitz-ian Mar-kov-ian Noe-ther-ian No-vem-ber Rie-mann-ian
Schwarz-schild Sep-tem-ber
form per-iods Uni-ver-si-ty cri-ti-sism for-ma-lism}
\Invalid@\nofrills
\Invalid@\usualspace
\newif\ifnofrills@
\def\nofrills@#1#2{\relaxnext@
  \DN@{\ifx\next\nofrills
    \nofrills@true\let#2\relax\DN@\nofrills{\nextii@}%
  \else
    \nofrills@false\def#2{#1}\let\next@\nextii@\fi
\next@}}
\def\usualspace@#1{\ifnofrills@\def\usualspace{#1}\fi}
\def\addto#1#2{\csname \expandafter\eat@\string#1@\endcsname
  \expandafter{\the\csname \expandafter\eat@\string#1@\endcsname#2}}
\newdimen\bigsize@
\def\big@#1#2{{\hbox{$\left#2\vcenter to#1\bigsize@{}%
  \right.\nulldelimiterspace\z@\m@th$}}}
\def\big{\big@\@ne}
\def\Big{\big@{1.5}}
\def\bigg{\big@\tw@}
\def\Bigg{\big@{2.5}}
\def\raggedcenter@{\leftskip\z@ plus.4\hsize \rightskip\leftskip
 \parfillskip\z@ \parindent\z@ \spaceskip.3333em \xspaceskip.5em
 \pretolerance9999\tolerance9999 \exhyphenpenalty\@M
 \hyphenpenalty\@M \let\\\linebreak}
\def\upperspecialchars{\def\ss{SS}\let\i=I\let\j=J\let\ae\AE\let\oe\OE
  \let\o\O\let\aa\AA\let\l\L}
\def\uppercasetext@#1{%
  {\spaceskip1.2\fontdimen2\the\font plus1.2\fontdimen3\the\font
   \upperspecialchars\uctext@#1$\m@th\aftergroup\eat@$}}
\def\uctext@#1$#2${\endash@#1-\endash@$#2$\uctext@}
\def\endash@#1-#2\endash@{\uppercase{#1}\if\notempty{#2}--\endash@#2\endash@\fi}
\def\runaway@#1{\DN@{#1}\ifx\envir@\next@
  \Err@{You seem to have a missing or misspelled \string\end#1 ...}%
  \let\envir@\empty\fi}
\newif\iftemp@
\def\notempty#1{TT\fi\def\test@{#1}\ifx\test@\empty\temp@false
  \else\temp@true\fi \iftemp@}
\font@\tensmc=cmcsc10
\font@\sevenex=cmex7
\font@\sevenit=cmti7
\font@\eightrm=cmr8 
\font@\sixrm=cmr6 
\font@\eighti=cmmi8     \skewchar\eighti='177 
\font@\sixi=cmmi6       \skewchar\sixi='177   
\font@\eightsy=cmsy8    \skewchar\eightsy='60 
\font@\sixsy=cmsy6      \skewchar\sixsy='60   
\font@\eightex=cmex8
\font@\eightbf=cmbx8 
\font@\sixbf=cmbx6   
\font@\eightit=cmti8 
\font@\eightsl=cmsl8 
\font@\eightsmc=cmcsc8
\font@\eighttt=cmtt8 


\loadmsam
\loadmsbm
\loadeufm
\UseAMSsymbols
\newtoks\tenpoint@
\def\tenpoint{\normalbaselineskip12\p@
 \abovedisplayskip12\p@ plus3\p@ minus9\p@
 \belowdisplayskip\abovedisplayskip
 \abovedisplayshortskip\z@ plus3\p@
 \belowdisplayshortskip7\p@ plus3\p@ minus4\p@
 \textonlyfont@\rm\tenrm \textonlyfont@\it\tenit
 \textonlyfont@\sl\tensl \textonlyfont@\bf\tenbf
 \textonlyfont@\smc\tensmc \textonlyfont@\tt\tentt
 \textonlyfont@\bsmc\tenbsmc
 \ifsyntax@ \def\big##1{{\hbox{$\left##1\right.$}}}%
  \let\Big\big \let\bigg\big \let\Bigg\big
 \else
  \textfont\z@=\tenrm  \scriptfont\z@=\sevenrm  \scriptscriptfont\z@=\fiverm
  \textfont\@ne=\teni  \scriptfont\@ne=\seveni  \scriptscriptfont\@ne=\fivei
  \textfont\tw@=\tensy \scriptfont\tw@=\sevensy \scriptscriptfont\tw@=\fivesy
  \textfont\thr@@=\tenex \scriptfont\thr@@=\sevenex
        \scriptscriptfont\thr@@=\sevenex
  \textfont\itfam=\tenit \scriptfont\itfam=\sevenit
        \scriptscriptfont\itfam=\sevenit
  \textfont\bffam=\tenbf \scriptfont\bffam=\sevenbf
        \scriptscriptfont\bffam=\fivebf
  \setbox\strutbox\hbox{\vrule height8.5\p@ depth3.5\p@ width\z@}%
  \setbox\strutbox@\hbox{\lower.5\normallineskiplimit\vbox{%
        \kern-\normallineskiplimit\copy\strutbox}}%
 \setbox\z@\vbox{\hbox{$($}\kern\z@}\bigsize@=1.2\ht\z@
 \fi
 \normalbaselines\rm\ex@.2326ex\jot3\ex@\the\tenpoint@}
\newtoks\eightpoint@
\def\eightpoint{\normalbaselineskip10\p@
 \abovedisplayskip10\p@ plus2.4\p@ minus7.2\p@
 \belowdisplayskip\abovedisplayskip
 \abovedisplayshortskip\z@ plus2.4\p@
 \belowdisplayshortskip5.6\p@ plus2.4\p@ minus3.2\p@
 \textonlyfont@\rm\eightrm \textonlyfont@\it\eightit
 \textonlyfont@\sl\eightsl \textonlyfont@\bf\eightbf
 \textonlyfont@\smc\eightsmc \textonlyfont@\tt\eighttt
 \textonlyfont@\bsmc\eightbsmc
 \ifsyntax@\def\big##1{{\hbox{$\left##1\right.$}}}%
  \let\Big\big \let\bigg\big \let\Bigg\big
 \else
  \textfont\z@=\eightrm \scriptfont\z@=\sixrm \scriptscriptfont\z@=\fiverm
  \textfont\@ne=\eighti \scriptfont\@ne=\sixi \scriptscriptfont\@ne=\fivei
  \textfont\tw@=\eightsy \scriptfont\tw@=\sixsy \scriptscriptfont\tw@=\fivesy
  \textfont\thr@@=\eightex \scriptfont\thr@@=\sevenex
   \scriptscriptfont\thr@@=\sevenex
  \textfont\itfam=\eightit \scriptfont\itfam=\sevenit
   \scriptscriptfont\itfam=\sevenit
  \textfont\bffam=\eightbf \scriptfont\bffam=\sixbf
   \scriptscriptfont\bffam=\fivebf
 \setbox\strutbox\hbox{\vrule height7\p@ depth3\p@ width\z@}%
 \setbox\strutbox@\hbox{\raise.5\normallineskiplimit\vbox{%
   \kern-\normallineskiplimit\copy\strutbox}}%
 \setbox\z@\vbox{\hbox{$($}\kern\z@}\bigsize@=1.2\ht\z@
 \fi
 \normalbaselines\eightrm\ex@.2326ex\jot3\ex@\the\eightpoint@}

\font@\twelverm=cmr10 scaled\magstep1
\font@\twelveit=cmti10 scaled\magstep1
\font@\twelvesl=cmsl10 scaled\magstep1
\font@\twelvesmc=cmcsc10 scaled\magstep1
\font@\twelvett=cmtt10 scaled\magstep1
\font@\twelvebf=cmbx10 scaled\magstep1
\font@\twelvei=cmmi10 scaled\magstep1
\font@\twelvesy=cmsy10 scaled\magstep1
\font@\twelveex=cmex10 scaled\magstep1
\font@\twelvemsa=msam10 scaled\magstep1
\font@\twelveeufm=eufm10 scaled\magstep1
\font@\twelvemsb=msbm10 scaled\magstep1
\newtoks\twelvepoint@
\def\twelvepoint{\normalbaselineskip15\p@
 \abovedisplayskip15\p@ plus3.6\p@ minus10.8\p@
 \belowdisplayskip\abovedisplayskip
 \abovedisplayshortskip\z@ plus3.6\p@
 \belowdisplayshortskip8.4\p@ plus3.6\p@ minus4.8\p@
 \textonlyfont@\rm\twelverm \textonlyfont@\it\twelveit
 \textonlyfont@\sl\twelvesl \textonlyfont@\bf\twelvebf
 \textonlyfont@\smc\twelvesmc \textonlyfont@\tt\twelvett
 \textonlyfont@\bsmc\twelvebsmc
 \ifsyntax@ \def\big##1{{\hbox{$\left##1\right.$}}}%
  \let\Big\big \let\bigg\big \let\Bigg\big
 \else
  \textfont\z@=\twelverm  \scriptfont\z@=\tenrm  \scriptscriptfont\z@=\sevenrm
  \textfont\@ne=\twelvei  \scriptfont\@ne=\teni  \scriptscriptfont\@ne=\seveni
  \textfont\tw@=\twelvesy \scriptfont\tw@=\tensy \scriptscriptfont\tw@=\sevensy
  \textfont\thr@@=\twelveex \scriptfont\thr@@=\tenex
        \scriptscriptfont\thr@@=\tenex
  \textfont\itfam=\twelveit \scriptfont\itfam=\tenit
        \scriptscriptfont\itfam=\tenit
  \textfont\bffam=\twelvebf \scriptfont\bffam=\tenbf
        \scriptscriptfont\bffam=\sevenbf
  \setbox\strutbox\hbox{\vrule height10.2\p@ depth4.2\p@ width\z@}%
  \setbox\strutbox@\hbox{\lower.6\normallineskiplimit\vbox{%
        \kern-\normallineskiplimit\copy\strutbox}}%
 \setbox\z@\vbox{\hbox{$($}\kern\z@}\bigsize@=1.4\ht\z@
 \fi
 \normalbaselines\rm\ex@.2326ex\jot3.6\ex@\the\twelvepoint@}

\def\headfonts{\twelvepoint\bf}

\font@\fourteenrm=cmr10 scaled\magstep2
\font@\fourteenit=cmti10 scaled\magstep2
\font@\fourteensl=cmsl10 scaled\magstep2
\font@\fourteensmc=cmcsc10 scaled\magstep2
\font@\fourteentt=cmtt10 scaled\magstep2
\font@\fourteenbf=cmbx10 scaled\magstep2
\font@\fourteeni=cmmi10 scaled\magstep2
\font@\fourteensy=cmsy10 scaled\magstep2
\font@\fourteenex=cmex10 scaled\magstep2
\font@\fourteenmsa=msam10 scaled\magstep2
\font@\fourteeneufm=eufm10 scaled\magstep2
\font@\fourteenmsb=msbm10 scaled\magstep2
\newtoks\fourteenpoint@
\def\fourteenpoint{\normalbaselineskip15\p@
 \abovedisplayskip18\p@ plus4.3\p@ minus12.9\p@
 \belowdisplayskip\abovedisplayskip
 \abovedisplayshortskip\z@ plus4.3\p@
 \belowdisplayshortskip10.1\p@ plus4.3\p@ minus5.8\p@
 \textonlyfont@\rm\fourteenrm \textonlyfont@\it\fourteenit
 \textonlyfont@\sl\fourteensl \textonlyfont@\bf\fourteenbf
 \textonlyfont@\smc\fourteensmc \textonlyfont@\tt\fourteentt
 \textonlyfont@\bsmc\fourteenbsmc
 \ifsyntax@ \def\big##1{{\hbox{$\left##1\right.$}}}%
  \let\Big\big \let\bigg\big \let\Bigg\big
 \else
  \textfont\z@=\fourteenrm  \scriptfont\z@=\twelverm  \scriptscriptfont\z@=\tenrm
  \textfont\@ne=\fourteeni  \scriptfont\@ne=\twelvei  \scriptscriptfont\@ne=\teni
  \textfont\tw@=\fourteensy \scriptfont\tw@=\twelvesy \scriptscriptfont\tw@=\tensy
  \textfont\thr@@=\fourteenex \scriptfont\thr@@=\twelveex
        \scriptscriptfont\thr@@=\twelveex
  \textfont\itfam=\fourteenit \scriptfont\itfam=\twelveit
        \scriptscriptfont\itfam=\twelveit
  \textfont\bffam=\fourteenbf \scriptfont\bffam=\twelvebf
        \scriptscriptfont\bffam=\tenbf
  \setbox\strutbox\hbox{\vrule height12.2\p@ depth5\p@ width\z@}%
  \setbox\strutbox@\hbox{\lower.72\normallineskiplimit\vbox{%
        \kern-\normallineskiplimit\copy\strutbox}}%
 \setbox\z@\vbox{\hbox{$($}\kern\z@}\bigsize@=1.7\ht\z@
 \fi
 \normalbaselines\rm\ex@.2326ex\jot4.3\ex@\the\fourteenpoint@}

\def\chapheadfonts{\fourteenpoint\bf}

\font@\seventeenrm=cmr10 scaled\magstep3
\font@\seventeenit=cmti10 scaled\magstep3
\font@\seventeensl=cmsl10 scaled\magstep3
\font@\seventeensmc=cmcsc10 scaled\magstep3
\font@\seventeentt=cmtt10 scaled\magstep3
\font@\seventeenbf=cmbx10 scaled\magstep3
\font@\seventeeni=cmmi10 scaled\magstep3
\font@\seventeensy=cmsy10 scaled\magstep3
\font@\seventeenex=cmex10 scaled\magstep3
\font@\seventeenmsa=msam10 scaled\magstep3
\font@\seventeeneufm=eufm10 scaled\magstep3
\font@\seventeenmsb=msbm10 scaled\magstep3
\newtoks\seventeenpoint@
\def\seventeenpoint{\normalbaselineskip18\p@
 \abovedisplayskip21.6\p@ plus5.2\p@ minus15.4\p@
 \belowdisplayskip\abovedisplayskip
 \abovedisplayshortskip\z@ plus5.2\p@
 \belowdisplayshortskip12.1\p@ plus5.2\p@ minus7\p@
 \textonlyfont@\rm\seventeenrm \textonlyfont@\it\seventeenit
 \textonlyfont@\sl\seventeensl \textonlyfont@\bf\seventeenbf
 \textonlyfont@\smc\seventeensmc \textonlyfont@\tt\seventeentt
 \textonlyfont@\bsmc\seventeenbsmc
 \ifsyntax@ \def\big##1{{\hbox{$\left##1\right.$}}}%
  \let\Big\big \let\bigg\big \let\Bigg\big
 \else
  \textfont\z@=\seventeenrm  \scriptfont\z@=\fourteenrm  \scriptscriptfont\z@=\twelverm
  \textfont\@ne=\seventeeni  \scriptfont\@ne=\fourteeni  \scriptscriptfont\@ne=\twelvei
  \textfont\tw@=\seventeensy \scriptfont\tw@=\fourteensy \scriptscriptfont\tw@=\twelvesy
  \textfont\thr@@=\seventeenex \scriptfont\thr@@=\fourteenex
        \scriptscriptfont\thr@@=\fourteenex
  \textfont\itfam=\seventeenit \scriptfont\itfam=\fourteenit
        \scriptscriptfont\itfam=\fourteenit
  \textfont\bffam=\seventeenbf \scriptfont\bffam=\fourteenbf
        \scriptscriptfont\bffam=\twelvebf
  \setbox\strutbox\hbox{\vrule height14.6\p@ depth6\p@ width\z@}%
  \setbox\strutbox@\hbox{\lower.86\normallineskiplimit\vbox{%
        \kern-\normallineskiplimit\copy\strutbox}}%
 \setbox\z@\vbox{\hbox{$($}\kern\z@}\bigsize@=2\ht\z@
 \fi
 \normalbaselines\rm\ex@.2326ex\jot5.2\ex@\the\seventeenpoint@}

\font@\rrrrrm=cmr10 scaled\magstep4
\font@\bigtitlefont=cmbx10 scaled\magstep4

\parindent1pc
\normallineskiplimit\p@
\newdimen\indenti \indenti=2pc
\def\pageheight#1{\vsize#1}
\def\pagewidth#1{\hsize#1%
   \captionwidth@\hsize \advance\captionwidth@-2\indenti}
\pagewidth{30pc} \pageheight{47pc}
\def\topmatter{%
 \ifx\undefined\msafam
 \else\font@\eightmsa=msam8 \font@\sixmsa=msam6
   \ifsyntax@\else \addto\tenpoint{\textfont\msafam=\tenmsa
              \scriptfont\msafam=\sevenmsa \scriptscriptfont\msafam=\fivemsa}%
     \addto\eightpoint{\textfont\msafam=\eightmsa \scriptfont\msafam=\sixmsa
              \scriptscriptfont\msafam=\fivemsa}%
   \fi
 \fi
 \ifx\undefined\msbfam
 \else\font@\eightmsb=msbm8 \font@\sixmsb=msbm6
   \ifsyntax@\else \addto\tenpoint{\textfont\msbfam=\tenmsb
         \scriptfont\msbfam=\sevenmsb \scriptscriptfont\msbfam=\fivemsb}%
     \addto\eightpoint{\textfont\msbfam=\eightmsb \scriptfont\msbfam=\sixmsb
         \scriptscriptfont\msbfam=\fivemsb}%
   \fi
 \fi
 \ifx\undefined\eufmfam
 \else \font@\eighteufm=eufm8 \font@\sixeufm=eufm6
   \ifsyntax@\else \addto\tenpoint{\textfont\eufmfam=\teneufm
       \scriptfont\eufmfam=\seveneufm \scriptscriptfont\eufmfam=\fiveeufm}%
     \addto\eightpoint{\textfont\eufmfam=\eighteufm
       \scriptfont\eufmfam=\sixeufm \scriptscriptfont\eufmfam=\fiveeufm}%
   \fi
 \fi
 \ifx\undefined\eufbfam
 \else \font@\eighteufb=eufb8 \font@\sixeufb=eufb6
   \ifsyntax@\else \addto\tenpoint{\textfont\eufbfam=\teneufb
      \scriptfont\eufbfam=\seveneufb \scriptscriptfont\eufbfam=\fiveeufb}%
    \addto\eightpoint{\textfont\eufbfam=\eighteufb
      \scriptfont\eufbfam=\sixeufb \scriptscriptfont\eufbfam=\fiveeufb}%
   \fi
 \fi
 \ifx\undefined\eusmfam
 \else \font@\eighteusm=eusm8 \font@\sixeusm=eusm6
   \ifsyntax@\else \addto\tenpoint{\textfont\eusmfam=\teneusm
       \scriptfont\eusmfam=\seveneusm \scriptscriptfont\eusmfam=\fiveeusm}%
     \addto\eightpoint{\textfont\eusmfam=\eighteusm
       \scriptfont\eusmfam=\sixeusm \scriptscriptfont\eusmfam=\fiveeusm}%
   \fi
 \fi
 \ifx\undefined\eusbfam
 \else \font@\eighteusb=eusb8 \font@\sixeusb=eusb6
   \ifsyntax@\else \addto\tenpoint{\textfont\eusbfam=\teneusb
       \scriptfont\eusbfam=\seveneusb \scriptscriptfont\eusbfam=\fiveeusb}%
     \addto\eightpoint{\textfont\eusbfam=\eighteusb
       \scriptfont\eusbfam=\sixeusb \scriptscriptfont\eusbfam=\fiveeusb}%
   \fi
 \fi
 \ifx\undefined\eurmfam
 \else \font@\eighteurm=eurm8 \font@\sixeurm=eurm6
   \ifsyntax@\else \addto\tenpoint{\textfont\eurmfam=\teneurm
       \scriptfont\eurmfam=\seveneurm \scriptscriptfont\eurmfam=\fiveeurm}%
     \addto\eightpoint{\textfont\eurmfam=\eighteurm
       \scriptfont\eurmfam=\sixeurm \scriptscriptfont\eurmfam=\fiveeurm}%
   \fi
 \fi
 \ifx\undefined\eurbfam
 \else \font@\eighteurb=eurb8 \font@\sixeurb=eurb6
   \ifsyntax@\else \addto\tenpoint{\textfont\eurbfam=\teneurb
       \scriptfont\eurbfam=\seveneurb \scriptscriptfont\eurbfam=\fiveeurb}%
    \addto\eightpoint{\textfont\eurbfam=\eighteurb
       \scriptfont\eurbfam=\sixeurb \scriptscriptfont\eurbfam=\fiveeurb}%
   \fi
 \fi
 \ifx\undefined\cmmibfam
 \else \font@\eightcmmib=cmmib8 \font@\sixcmmib=cmmib6
   \ifsyntax@\else \addto\tenpoint{\textfont\cmmibfam=\tencmmib
       \scriptfont\cmmibfam=\sevencmmib \scriptscriptfont\cmmibfam=\fivecmmib}%
    \addto\eightpoint{\textfont\cmmibfam=\eightcmmib
       \scriptfont\cmmibfam=\sixcmmib \scriptscriptfont\cmmibfam=\fivecmmib}%
   \fi
 \fi
 \ifx\undefined\cmbsyfam
 \else \font@\eightcmbsy=cmbsy8 \font@\sixcmbsy=cmbsy6
   \ifsyntax@\else \addto\tenpoint{\textfont\cmbsyfam=\tencmbsy
      \scriptfont\cmbsyfam=\sevencmbsy \scriptscriptfont\cmbsyfam=\fivecmbsy}%
    \addto\eightpoint{\textfont\cmbsyfam=\eightcmbsy
      \scriptfont\cmbsyfam=\sixcmbsy \scriptscriptfont\cmbsyfam=\fivecmbsy}%
   \fi
 \fi
 \let\topmatter\relax}
\def\chapterno@{\uppercase\expandafter{\romannumeral\chaptercount@}}
\newcount\chaptercount@
\def\chapter{\nofrills@{\afterassignment\chapterno@
                        CHAPTER \global\chaptercount@=}\chapter@
 \DNii@##1{\leavevmode\hskip-\leftskip
   \rlap{\vbox to\z@{\vss\centerline{\eightpoint
   \chapter@##1\unskip}\baselineskip2pc\null}}\hskip\leftskip
   \nofrills@false}%
 \FN@\next@}
\newbox\titlebox@

\def\title{\nofrills@{\relax}\title@%
 \DNii@##1\endtitle{\global\setbox\titlebox@\vtop{\tenpoint\bf
 \raggedcenter@\ignorespaces
 \baselineskip1.3\baselineskip\title@{##1}\endgraf}%
 \ifmonograph@ \edef\next{\the\leftheadtoks}\ifx\next\empty
    \leftheadtext{##1}\fi
 \fi
 \edef\next{\the\rightheadtoks}\ifx\next\empty \rightheadtext{##1}\fi
 }\FN@\next@}
\newbox\authorbox@
\def\author#1\endauthor{\global\setbox\authorbox@
 \vbox{\tenpoint\smc\raggedcenter@\ignorespaces
 #1\endgraf}\relaxnext@ \edef\next{\the\leftheadtoks}%
 \ifx\next\empty\leftheadtext{#1}\fi}
\newbox\affilbox@
\def\affil#1\endaffil{\global\setbox\affilbox@
 \vbox{\tenpoint\raggedcenter@\ignorespaces#1\endgraf}}
\newcount\addresscount@
\addresscount@\z@
\def\address#1\endaddress{\global\advance\addresscount@\@ne
  \expandafter\gdef\csname address\number\addresscount@\endcsname
  {\vskip12\p@ minus6\p@\noindent\eightpoint\smc\ignorespaces#1\par}}
\def\email{\nofrills@{\eightpoint{\it E-mail\/}:\enspace}\email@
  \DNii@##1\endemail{%
  \expandafter\gdef\csname email\number\addresscount@\endcsname
  {\def\usualspace{{\it\enspace}}\smallskip\noindent\eightpoint\email@
  \ignorespaces##1\par}}%
 \FN@\next@}
\def\thedate@{}
\def\date#1\enddate{\gdef\thedate@{\tenpoint\ignorespaces#1\unskip}}
\def\thethanks@{}
\def\thanks#1\endthanks{\gdef\thethanks@{\eightpoint\ignorespaces#1.\unskip}}
\def\thekeywords@{}
\def\keywords{\nofrills@{{\it Key words and phrases.\enspace}}\keywords@
 \DNii@##1\endkeywords{\def\thekeywords@{\def\usualspace{{\it\enspace}}%
 \eightpoint\keywords@\ignorespaces##1\unskip.}}%
 \FN@\next@}
\def\thesubjclass@{}
\def\subjclass{\nofrills@{{\rm2010 {\it Mathematics Subject
   Classification\/}.\enspace}}\subjclass@
 \DNii@##1\endsubjclass{\def\thesubjclass@{\def\usualspace
  {{\rm\enspace}}\eightpoint\subjclass@\ignorespaces##1\unskip.}}%
 \FN@\next@}
\newbox\abstractbox@
\def\abstract{\nofrills@{{\smc Abstract.\enspace}}\abstract@
 \DNii@{\setbox\abstractbox@\vbox\bgroup\noindent$$\vbox\bgroup
  \def\envir@{abstract}\advance\hsize-2\indenti
  \usualspace@{{\enspace}}\eightpoint \noindent\abstract@\ignorespaces}%
 \FN@\next@}
\def\endabstract{\par\unskip\egroup$$\egroup}
\def\widestnumber#1#2{\begingroup\let\head\null\let\subhead\empty
   \let\subsubhead\subhead
   \ifx#1\head\global\setbox\tocheadbox@\hbox{#2.\enspace}%
   \else\ifx#1\subhead\global\setbox\tocsubheadbox@\hbox{#2.\enspace}%
   \else\ifx#1\key\bgroup\let\endrefitem@\egroup
        \key#2\endrefitem@\global\refindentwd\wd\keybox@
   \else\ifx#1\no\bgroup\let\endrefitem@\egroup
        \no#2\endrefitem@\global\refindentwd\wd\nobox@
   \else\ifx#1\page\global\setbox\pagesbox@\hbox{\quad\bf#2}%
   \else\ifx#1\item\setboxz@h{#2}\global\rosteritemwd\wdz@
        \global\advance\rosteritemwd by.5\parindent
   \else\message{\string\widestnumber is not defined for this option
   (\string#1)}%
\fi\fi\fi\fi\fi\fi\endgroup}
\newif\ifmonograph@
\def\Monograph{\monograph@true \let\headmark\rightheadtext
  \let\varindent@\indent \def\headfont@{\bf}\def\proclaimheadfont@{\smc}%
  \def\demofont@{\smc}}
\let\varindent@\indent

\newbox\tocheadbox@    \newbox\tocsubheadbox@
\newbox\tocbox@
\def\toc{\toc@{Contents}}
\def\newtocdefs{%
   \def \title##1\endtitle
       {\penaltyandskip@\z@\smallskipamount
        \hangindent\wd\tocheadbox@\noindent{\bf##1}}%
   \def \chapter##1{%
        Chapter \uppercase\expandafter{\romannumeral##1.\unskip}\enspace}%
   \def \specialhead##1\endspecialhead
       {\par\hangindent\wd\tocheadbox@ \noindent##1\par}%
   \def \head##1 ##2\endhead
       {\par\hangindent\wd\tocheadbox@ \noindent
        \if\notempty{##1}\hbox to\wd\tocheadbox@{\hfil##1\enspace}\fi
        ##2\par}%
   \def \subhead##1 ##2\endsubhead
       {\par\vskip-\parskip {\normalbaselines
        \advance\leftskip\wd\tocheadbox@
        \hangindent\wd\tocsubheadbox@ \noindent
        \if\notempty{##1}\hbox to\wd\tocsubheadbox@{##1\unskip\hfil}\fi
         ##2\par}}%
   \def \subsubhead##1 ##2\endsubsubhead
       {\par\vskip-\parskip {\normalbaselines
        \advance\leftskip\wd\tocheadbox@
        \hangindent\wd\tocsubheadbox@ \noindent
        \if\notempty{##1}\hbox to\wd\tocsubheadbox@{##1\unskip\hfil}\fi
        ##2\par}}}
\def\toc@#1{\relaxnext@
   \def\page##1%
       {\unskip\penalty0\null\hfil
        \rlap{\hbox to\wd\pagesbox@{\quad\hfil##1}}\hfilneg\penalty\@M}%
 \DN@{\ifx\next\nofrills\DN@\nofrills{\nextii@}%
      \else\DN@{\nextii@{{#1}}}\fi
      \next@}%
 \DNii@##1{%
\ifmonograph@\bgroup\else\setbox\tocbox@\vbox\bgroup
   \centerline{\headfont@\ignorespaces##1\unskip}\nobreak
   \vskip\belowheadskip \fi
   \setbox\tocheadbox@\hbox{0.\enspace}%
   \setbox\tocsubheadbox@\hbox{0.0.\enspace}%
   \leftskip\indenti \rightskip\leftskip
   \setbox\pagesbox@\hbox{\bf\quad000}\advance\rightskip\wd\pagesbox@
   \newtocdefs
 }%
 \FN@\next@}
\def\endtoc{\par\egroup}
\let\pretitle\relax
\let\preauthor\relax
\let\preaffil\relax
\let\predate\relax
\let\preabstract\relax
\let\prepaper\relax
\def\dedicatory #1\enddedicatory{\def\preabstract{{\medskip
  \eightpoint\it \raggedcenter@#1\endgraf}}}
\def\thetranslator@{}
\def\translator#1\endtranslator{\def\thetranslator@{\nobreak\medskip
 \line{\eightpoint\hfil Translated by \uppercase{#1}\qquad\qquad}\nobreak}}
\outer\def\endtopmatter{\runaway@{abstract}%
 \edef\next{\the\leftheadtoks}\ifx\next\empty
  \expandafter\leftheadtext\expandafter{\the\rightheadtoks}\fi
 \ifmonograph@\else
   \ifx\thesubjclass@\empty\else \makefootnote@{}{\thesubjclass@}\fi
   \ifx\thekeywords@\empty\else \makefootnote@{}{\thekeywords@}\fi
   \ifx\thethanks@\empty\else \makefootnote@{}{\thethanks@}\fi
 \fi
  \pretitle
  \ifmonograph@ \topskip7pc \else \topskip4pc \fi
  \box\titlebox@
  \topskip10pt
  \preauthor
  \ifvoid\authorbox@\else \vskip2.5pc plus1pc \unvbox\authorbox@\fi
  \preaffil
  \ifvoid\affilbox@\else \vskip1pc plus.5pc \unvbox\affilbox@\fi
  \predate
  \ifx\thedate@\empty\else \vskip1pc plus.5pc \line{\hfil\thedate@\hfil}\fi
  \preabstract
  \ifvoid\abstractbox@\else \vskip1.5pc plus.5pc \unvbox\abstractbox@ \fi
  \ifvoid\tocbox@\else\vskip1.5pc plus.5pc \unvbox\tocbox@\fi
  \prepaper
  \vskip2pc plus1pc
}
\def\document{\let\fontlist@\relax\let\alloclist@\relax
  \tenpoint}

\newskip\aboveheadskip       \aboveheadskip1.8\bigskipamount
\newdimen\belowheadskip      \belowheadskip1.8\medskipamount

\def\headfont@{\smc}
\def\penaltyandskip@#1#2{\relax\ifdim\lastskip<#2\relax\removelastskip
      \ifnum#1=\z@\else\penalty@#1\relax\fi\vskip#2%
  \else\ifnum#1=\z@\else\penalty@#1\relax\fi\fi}
\def\nobreak{\penalty\@M
  \ifvmode\def\penalty@{\let\penalty@\penalty\count@@@}%
  \everypar{\let\penalty@\penalty\everypar{}}\fi}
\let\penalty@\penalty
\def\heading#1\endheading{\head#1\endhead}

\def\specialheadfont@{\bf}
\outer\def\specialhead{\par\penaltyandskip@{-200}\aboveheadskip
  \begingroup\interlinepenalty\@M\rightskip\z@ plus\hsize \let\\\linebreak
  \specialheadfont@\noindent\ignorespaces}
\def\endspecialhead{\par\endgroup\nobreak\vskip\belowheadskip}
\let\headmark\eat@
\newskip\subheadskip       \subheadskip\medskipamount
\def\subheadfont@{\bf}
\outer\def\subhead{\nofrills@{.\enspace}\subhead@
 \DNii@##1\endsubhead{\par\penaltyandskip@{-100}\subheadskip
  \varindent@{\usualspace@{{\subheadfont@\enspace}}%
 \subheadfont@\ignorespaces##1\unskip\subhead@}\ignorespaces}%
 \FN@\next@}
\outer\def\subsubhead{\nofrills@{.\enspace}\subsubhead@
 \DNii@##1\endsubsubhead{\par\penaltyandskip@{-50}\medskipamount
      {\usualspace@{{\it\enspace}}%
  \it\ignorespaces##1\unskip\subsubhead@}\ignorespaces}%
 \FN@\next@}
\def\proclaimheadfont@{\bf}
\outer\def\proclaim{\runaway@{proclaim}\def\envir@{proclaim}%
  \nofrills@{.\enspace}\proclaim@
 \DNii@##1{\penaltyandskip@{-100}\medskipamount\varindent@
   \usualspace@{{\proclaimheadfont@\enspace}}\proclaimheadfont@
   \ignorespaces##1\unskip\proclaim@
  \sl\ignorespaces}%
 \FN@\next@}
\outer\def\endproclaim{\let\envir@\relax\par\rm
  \penaltyandskip@{55}\medskipamount}
\def\demoheadfont@{\it}
\def\demo{\runaway@{proclaim}\nofrills@{.\enspace}\demo@
     \DNii@##1{\par\penaltyandskip@\z@\medskipamount
  {\usualspace@{{\demoheadfont@\enspace}}%
  \varindent@\demoheadfont@\ignorespaces##1\unskip\demo@}\rm
  \ignorespaces}\FN@\next@}
\def\enddemo{\par\medskip}
\def\qed{\ifhmode\unskip\nobreak\fi\quad\ifmmode\square\else$\m@th\square$\fi}
\let\remark\demo
\let\endremark\enddemo
\def\definition{\runaway@{proclaim}%
  \nofrills@{.\demoheadfont@\enspace}\definition@
        \DNii@##1{\penaltyandskip@{-100}\medskipamount
        {\usualspace@{{\demoheadfont@\enspace}}%
        \varindent@\demoheadfont@\ignorespaces##1\unskip\definition@}%
        \rm \ignorespaces}\FN@\next@}


\newdimen\rosteritemwd
\newcount\rostercount@
\newif\iffirstitem@
\let\plainitem@\item
\newtoks\everypartoks@
\def\par@{\everypartoks@\expandafter{\the\everypar}\everypar{}}
\def\roster{\edef\leftskip@{\leftskip\the\leftskip}%
 \relaxnext@
 \rostercount@\z@  
 \def\item{\FN@\rosteritem@}%
 \DN@{\ifx\next\runinitem\let\next@\nextii@\else
  \let\next@\nextiii@\fi\next@}%
 \DNii@\runinitem  
  {\unskip  
   \DN@{\ifx\next[\let\next@\nextii@\else
    \ifx\next"\let\next@\nextiii@\else\let\next@\nextiv@\fi\fi\next@}%
   \DNii@[####1]{\rostercount@####1\relax
    \enspace{\rm(\number\rostercount@)}~\ignorespaces}%
   \def\nextiii@"####1"{\enspace{\rm####1}~\ignorespaces}%
   \def\nextiv@{\enspace{\rm(1)}\rostercount@\@ne~}%
   \par@\firstitem@false  
   \FN@\next@}%
 \def\nextiii@{\par\par@  
  \penalty\@m\smallskip\vskip-\parskip
  \firstitem@true}%
 \FN@\next@}
\def\rosteritem@{\iffirstitem@\firstitem@false\else\par\vskip-\parskip\fi
 \leftskip3\parindent\noindent  
 \DNii@[##1]{\rostercount@##1\relax
  \llap{\hbox to2.5\parindent{\hss\rm(\number\rostercount@)}%
   \hskip.5\parindent}\ignorespaces}%
 \def\nextiii@"##1"{%
  \llap{\hbox to2.5\parindent{\hss\rm##1}\hskip.5\parindent}\ignorespaces}%
 \def\nextiv@{\advance\rostercount@\@ne
  \llap{\hbox to2.5\parindent{\hss\rm(\number\rostercount@)}%
   \hskip.5\parindent}}%
 \ifx\next[\let\next@\nextii@\else\ifx\next"\let\next@\nextiii@\else
  \let\next@\nextiv@\fi\fi\next@}

\newif\ifnextRunin@
\def\endroster{\relaxnext@
 \par\leftskip@  
 \penalty-50 \vskip-\parskip\smallskip  
 \DN@{\ifx\next\Runinitem\let\next@\relax
  \else\nextRunin@false\let\item\plainitem@  
   \ifx\next\par 
    \DN@\par{\everypar\expandafter{\the\everypartoks@}}%
   \else  
    \DN@{\noindent\everypar\expandafter{\the\everypartoks@}}%
  \fi\fi\next@}%
 \FN@\next@}
\newcount\rosterhangafter@
\def\Runinitem#1\roster\runinitem{\relaxnext@
 \rostercount@\z@ 
 \def\item{\FN@\rosteritem@}%
 \def\runinitem@{#1}%
 \DN@{\ifx\next[\let\next\nextii@\else\ifx\next"\let\next\nextiii@
  \else\let\next\nextiv@\fi\fi\next}%
 \DNii@[##1]{\rostercount@##1\relax
  \def\item@{{\rm(\number\rostercount@)}}\nextv@}%
 \def\nextiii@"##1"{\def\item@{{\rm##1}}\nextv@}%
 \def\nextiv@{\advance\rostercount@\@ne
  \def\item@{{\rm(\number\rostercount@)}}\nextv@}%
 \def\nextv@{\setbox\z@\vbox  
  {\ifnextRunin@\noindent\fi  
  \runinitem@\unskip\enspace\item@~\par  
  \global\rosterhangafter@\prevgraf}%
  \firstitem@false  
  \ifnextRunin@\else\par\fi  
  \hangafter\rosterhangafter@\hangindent3\parindent
  \ifnextRunin@\noindent\fi  
  \runinitem@\unskip\enspace 
  \item@~\ifnextRunin@\else\par@\fi  
  \nextRunin@true\ignorespaces}%
 \FN@\next@}
\def\footmarkform@#1{$\m@th^{#1}$}
\let\thefootnotemark\footmarkform@
\def\makefootnote@#1#2{\insert\footins
 {\interlinepenalty\interfootnotelinepenalty
 \eightpoint\splittopskip\ht\strutbox\splitmaxdepth\dp\strutbox
 \floatingpenalty\@MM\leftskip\z@\rightskip\z@\spaceskip\z@\xspaceskip\z@
 \leavevmode{#1}\footstrut\ignorespaces#2\unskip\lower\dp\strutbox
 \vbox to\dp\strutbox{}}}
\newcount\footmarkcount@
\footmarkcount@\z@
\def\footnotemark{\let\@sf\empty\relaxnext@
 \ifhmode\edef\@sf{\spacefactor\the\spacefactor}\/\fi
 \DN@{\ifx[\next\let\next@\nextii@\else
  \ifx"\next\let\next@\nextiii@\else
  \let\next@\nextiv@\fi\fi\next@}%
 \DNii@[##1]{\footmarkform@{##1}\@sf}%
 \def\nextiii@"##1"{{##1}\@sf}%
 \def\nextiv@{\iffirstchoice@\global\advance\footmarkcount@\@ne\fi
  \footmarkform@{\number\footmarkcount@}\@sf}%
 \FN@\next@}
\def\footnotetext{\relaxnext@
 \DN@{\ifx[\next\let\next@\nextii@\else
  \ifx"\next\let\next@\nextiii@\else
  \let\next@\nextiv@\fi\fi\next@}%
 \DNii@[##1]##2{\makefootnote@{\footmarkform@{##1}}{##2}}%
 \def\nextiii@"##1"##2{\makefootnote@{##1}{##2}}%
 \def\nextiv@##1{\makefootnote@{\footmarkform@{\number\footmarkcount@}}{##1}}%
 \FN@\next@}
\def\footnote{\let\@sf\empty\relaxnext@
 \ifhmode\edef\@sf{\spacefactor\the\spacefactor}\/\fi
 \DN@{\ifx[\next\let\next@\nextii@\else
  \ifx"\next\let\next@\nextiii@\else
  \let\next@\nextiv@\fi\fi\next@}%
 \DNii@[##1]##2{\footnotemark[##1]\footnotetext[##1]{##2}}%
 \def\nextiii@"##1"##2{\footnotemark"##1"\footnotetext"##1"{##2}}%
 \def\nextiv@##1{\footnotemark\footnotetext{##1}}%
 \FN@\next@}
\def\adjustfootnotemark#1{\advance\footmarkcount@#1\relax}
\def\footnoterule{\kern-3\p@
  \hrule width 5pc\kern 2.6\p@} 
\def\captionfont@{\smc}
\def\topcaption#1#2\endcaption{%
  {\dimen@\hsize \advance\dimen@-\captionwidth@
   \rm\raggedcenter@ \advance\leftskip.5\dimen@ \rightskip\leftskip
  {\captionfont@#1}%
  \if\notempty{#2}.\enspace\ignorespaces#2\fi
  \endgraf}\nobreak\bigskip}
\def\botcaption#1#2\endcaption{%
  \nobreak\bigskip
  \setboxz@h{\captionfont@#1\if\notempty{#2}.\enspace\rm#2\fi}%
  {\dimen@\hsize \advance\dimen@-\captionwidth@
   \leftskip.5\dimen@ \rightskip\leftskip
   \noindent \ifdim\wdz@>\captionwidth@ 
   \else\hfil\fi 
  {\captionfont@#1}\if\notempty{#2}.\enspace\rm#2\fi\endgraf}}
\def\@ins{\par\begingroup\def\vspace##1{\vskip##1\relax}%
  \def\captionwidth##1{\captionwidth@##1\relax}%
  \setbox\z@\vbox\bgroup} 
\def\block{\RIfMIfI@\nondmatherr@\block\fi
       \else\ifvmode\vskip\abovedisplayskip\noindent\fi
        $$\def\endblock{\par\egroup$$}\fi
  \vbox\bgroup\advance\hsize-2\indenti\noindent}
\def\endblock{\par\egroup}
\def\cite#1{{\rm[{\citefont@\m@th#1}]}}
\def\citefont@{\rm}
\def\refsfont@{\eightpoint}
\outer\def\Refs{\runaway@{proclaim}%
 \relaxnext@ \DN@{\ifx\next\nofrills\DN@\nofrills{\nextii@}\else
  \DN@{\nextii@{References}}\fi\next@}%
 \DNii@##1{\penaltyandskip@{-200}\aboveheadskip
  \line{\hfil\headfont@\ignorespaces##1\unskip\hfil}\nobreak
  \vskip\belowheadskip
  \begingroup\refsfont@\sfcode`.=\@m}%
 \FN@\next@}
\def\endRefs{\par\endgroup}
\newbox\nobox@            \newbox\keybox@           \newbox\bybox@
\newbox\paperbox@         \newbox\paperinfobox@     \newbox\jourbox@
\newbox\volbox@           \newbox\issuebox@         \newbox\yrbox@
\newbox\pagesbox@         \newbox\bookbox@          \newbox\bookinfobox@
\newbox\publbox@          \newbox\publaddrbox@      \newbox\finalinfobox@
\newbox\edsbox@           \newbox\langbox@
\newif\iffirstref@        \newif\iflastref@
\newif\ifprevjour@        \newif\ifbook@            \newif\ifprevinbook@
\newif\ifquotes@          \newif\ifbookquotes@      \newif\ifpaperquotes@
\newdimen\bysamerulewd@
\setboxz@h{\refsfont@\kern3em}
\bysamerulewd@\wdz@
\newdimen\refindentwd
\setboxz@h{\refsfont@ 00. }
\refindentwd\wdz@
\outer\def\ref{\begingroup \noindent\hangindent\refindentwd
 \firstref@true \def\nofrills{\def\refkern@{\kern3sp}}%
 \ref@}
\def\ref@{\book@false \bgroup\let\endrefitem@\egroup \ignorespaces}
\def\moreref{\endrefitem@\endref@\firstref@false\ref@}%
\def\transl{\endrefitem@\endref@\firstref@false
  \book@false
  \prepunct@
  \setboxz@h\bgroup \aftergroup\unhbox\aftergroup\z@
    \def\endrefitem@{\unskip\refkern@\egroup}\ignorespaces}%
\def\emptyifempty@{\dimen@\wd\currbox@
  \advance\dimen@-\wd\z@ \advance\dimen@-.1\p@
  \ifdim\dimen@<\z@ \setbox\currbox@\copy\voidb@x \fi}
\let\refkern@\relax
\def\endrefitem@{\unskip\refkern@\egroup
  \setboxz@h{\refkern@}\emptyifempty@}\ignorespaces
\def\refdef@#1#2#3{\edef\next@{\noexpand\endrefitem@
  \let\noexpand\currbox@\csname\expandafter\eat@\string#1box@\endcsname
    \noexpand\setbox\noexpand\currbox@\hbox\bgroup}%
  \toks@\expandafter{\next@}%
  \if\notempty{#2#3}\toks@\expandafter{\the\toks@
  \def\endrefitem@{\unskip#3\refkern@\egroup
  \setboxz@h{#2#3\refkern@}\emptyifempty@}#2}\fi
  \toks@\expandafter{\the\toks@\ignorespaces}%
  \edef#1{\the\toks@}}
\refdef@\no{}{. }
\refdef@\key{[\m@th}{] }
\refdef@\by{}{}
\def\bysame{\by\hbox to\bysamerulewd@{\hrulefill}\thinspace
   \kern0sp}
\def\manyby{\message{\string\manyby is no longer necessary; \string\by
  can be used instead, starting with version 2.0 of \styname.STY}\by}
\refdef@\paper{\ifpaperquotes@``\fi\it}{}
\refdef@\paperinfo{}{}
\def\jour{\endrefitem@\let\currbox@\jourbox@
  \setbox\currbox@\hbox\bgroup
  \def\endrefitem@{\unskip\refkern@\egroup
    \setboxz@h{\refkern@}\emptyifempty@
    \ifvoid\jourbox@\else\prevjour@true\fi}%
\ignorespaces}
\refdef@\vol{\ifbook@\else\bf\fi}{}
\refdef@\issue{no. }{}
\refdef@\yr{}{}
\refdef@\pages{}{}
\def\page{\endrefitem@\def\pp@{\def\pp@{pp.~}p.~}\let\currbox@\pagesbox@
  \setbox\currbox@\hbox\bgroup\ignorespaces}
\def\pp@{pp.~}
\def\book{\endrefitem@ \let\currbox@\bookbox@
 \setbox\currbox@\hbox\bgroup\def\endrefitem@{\unskip\refkern@\egroup
  \setboxz@h{\ifbookquotes@``\fi}\emptyifempty@
  \ifvoid\bookbox@\else\book@true\fi}%
  \ifbookquotes@``\fi\it\ignorespaces}
\def\inbook{\endrefitem@
  \let\currbox@\bookbox@\setbox\currbox@\hbox\bgroup
  \def\endrefitem@{\unskip\refkern@\egroup
  \setboxz@h{\ifbookquotes@``\fi}\emptyifempty@
  \ifvoid\bookbox@\else\book@true\previnbook@true\fi}%
  \ifbookquotes@``\fi\ignorespaces}
\refdef@\eds{(}{, eds.)}
\def\ed{\endrefitem@\let\currbox@\edsbox@
 \setbox\currbox@\hbox\bgroup
 \def\endrefitem@{\unskip, ed.)\refkern@\egroup
  \setboxz@h{(, ed.)}\emptyifempty@}(\ignorespaces}
\refdef@\bookinfo{}{}
\refdef@\publ{}{}
\refdef@\publaddr{}{}
\refdef@\finalinfo{}{}
\refdef@\lang{(}{)}

\let\refdef@\relax 
\def\ppunbox@#1{\ifvoid#1\else\prepunct@\unhbox#1\fi}
\def\nocomma@#1{\ifvoid#1\else\changepunct@3\prepunct@\unhbox#1\fi}
\def\changepunct@#1{\ifnum\lastkern<3 \unkern\kern#1sp\fi}
\def\prepunct@{\count@\lastkern\unkern
  \ifnum\lastpenalty=0
    \let\penalty@\relax
  \else
    \edef\penalty@{\penalty\the\lastpenalty\relax}%
  \fi
  \unpenalty
  \let\refspace@\ \ifcase\count@,
\or;\or.\or 
  \or\let\refspace@\relax
  \else,\fi
  \ifquotes@''\quotes@false\fi \penalty@ \refspace@
}
\def\transferpenalty@#1{\dimen@\lastkern\unkern
  \ifnum\lastpenalty=0\unpenalty\let\penalty@\relax
  \else\edef\penalty@{\penalty\the\lastpenalty\relax}\unpenalty\fi
  #1\penalty@\kern\dimen@}
\def\endref{\endrefitem@\lastref@true\endref@
  \par\endgroup \prevjour@false \previnbook@false }
\def\endref@{%
\iffirstref@
  \ifvoid\nobox@\ifvoid\keybox@\indent\fi
  \else\hbox to\refindentwd{\hss\unhbox\nobox@}\fi
  \ifvoid\keybox@
  \else\ifdim\wd\keybox@>\refindentwd
         \box\keybox@
       \else\hbox to\refindentwd{\unhbox\keybox@\hfil}\fi\fi
  \kern4sp\ppunbox@\bybox@
\fi 
  \ifvoid\paperbox@
  \else\prepunct@\unhbox\paperbox@
    \ifpaperquotes@\quotes@true\fi\fi
  \ppunbox@\paperinfobox@
  \ifvoid\jourbox@
    \ifprevjour@ \nocomma@\volbox@
      \nocomma@\issuebox@
      \ifvoid\yrbox@\else\changepunct@3\prepunct@(\unhbox\yrbox@
        \transferpenalty@)\fi
      \ppunbox@\pagesbox@
    \fi 
  \else \prepunct@\unhbox\jourbox@
    \nocomma@\volbox@
    \nocomma@\issuebox@
    \ifvoid\yrbox@\else\changepunct@3\prepunct@(\unhbox\yrbox@
      \transferpenalty@)\fi
    \ppunbox@\pagesbox@
  \fi 
  \ifbook@\prepunct@\unhbox\bookbox@ \ifbookquotes@\quotes@true\fi \fi
  \nocomma@\edsbox@
  \ppunbox@\bookinfobox@
  \ifbook@\ifvoid\volbox@\else\prepunct@ vol.~\unhbox\volbox@
  \fi\fi
  \ppunbox@\publbox@ \ppunbox@\publaddrbox@
  \ifbook@ \ppunbox@\yrbox@
    \ifvoid\pagesbox@
    \else\prepunct@\pp@\unhbox\pagesbox@\fi
  \else
    \ifprevinbook@ \ppunbox@\yrbox@
      \ifvoid\pagesbox@\else\prepunct@\pp@\unhbox\pagesbox@\fi
    \fi \fi
  \ppunbox@\finalinfobox@
  \iflastref@
    \ifvoid\langbox@.\ifquotes@''\fi
    \else\changepunct@2\prepunct@\unhbox\langbox@\fi
  \else
    \ifvoid\langbox@\changepunct@1%
    \else\changepunct@3\prepunct@\unhbox\langbox@
      \changepunct@1\fi
  \fi
}
\outer\def\enddocument{%
 \runaway@{proclaim}%
\ifmonograph@ 
\else
 \nobreak
 \thetranslator@
 \count@\z@ \loop\ifnum\count@<\addresscount@\advance\count@\@ne
 \csname address\number\count@\endcsname
 \csname email\number\count@\endcsname
 \repeat
\fi
 \vfill\supereject\end}

\def\headfont@{\headfonts}
\def\proclaimheadfont@{\bf}
\def\specialheadfont@{\bf}
\def\subheadfont@{\bf}
\def\demoheadfont@{\smc}

\newif\ifThisToToc \ThisToTocfalse
\newif\iftocloaded \tocloadedfalse

\def\C@L{\noexpand\Cal}\def\B@B{\noexpand\Bbb}\def\fR@K{\noexpand\frak}
\def\S@{\noexpand\S}\def\P@P{\noexpand\"}
\def\xpar{\\}

\def\writetoc#1{\iftocloaded\ifThisToToc\begingroup\def\totoc{}
  \def\Cal{\noexpand\C@L}\def\Bbb{\noexpand\B@B}
  \def\frak{\noexpand\fR@K}\def\goth{\frak}\def\S{\noexpand\S@}
  \def\"{\noexpand\P@P}
  \def\xpar{\par\penalty100000 }\def\idx##1{##1}\def\\{\xpar}
  \edef\next@{\write\toc{\noindent#1\leaderfill\noexpand\folio\par}}%
  \next@\endgroup\global\ThisToTocfalse\fi\fi}
\def\leaderfill{\leaders\hbox to 1em{\hss.\hss}\hfill}

\newif\ifindexloaded \indexloadedfalse
\def\idx#1{\ifindexloaded\begingroup\def\ign{}\def\it{}\def\/{}%
 \def\smc{}\def\bf{}\def\tt{}%
 \def\Cal{\noexpand\C@L}\def\Bbb{\noexpand\B@B}%
 \def\frak{\noexpand\fR@K}\def\goth{\frak}\def\S{\noexpand\S@}%
  \def\"{\noexpand\P@P}%
 {\edef\next@{\write\index{#1, \noexpand\folio}}\next@}%
 \endgroup\fi{#1}}
\def\ign#1{}

\def\totoc{\global\ThisToToctrue}

\outer\def\head#1\endhead{\par\penaltyandskip@{-200}\aboveheadskip
 {\headfont@\raggedcenter@\interlinepenalty\@M
 \ignorespaces#1\endgraf}\nobreak
 \vskip\belowheadskip
 \headmark{#1}\writetoc{#1}}

\outer\def\chaphead#1\endchaphead{\par\penaltyandskip@{-200}\aboveheadskip
 {\chapheadfonts\raggedcenter@\interlinepenalty\@M
 \ignorespaces#1\endgraf}\nobreak
 \vskip3\belowheadskip
 \headmark{#1}\writetoc{#1}}

\def\folio{{\foliofont@\ifnum\pageno<\z@ \romannumeral-\pageno
 \else\number\pageno \fi}}
\newtoks\leftheadtoks
\newtoks\rightheadtoks

\def\leftheadtext{\nofrills@{\relax}\lht@
  \DNii@##1{\leftheadtoks\expandafter{\lht@{##1}}%
    \mark{\the\leftheadtoks\noexpand\else\the\rightheadtoks}
    \ifsyntax@\setboxz@h{\def\\{\unskip\space\ignorespaces}%
        \headlinefont@##1}\fi}%
  \FN@\next@}
\def\rightheadtext{\nofrills@{\relax}\rht@
  \DNii@##1{\rightheadtoks\expandafter{\rht@{##1}}%
    \mark{\the\leftheadtoks\noexpand\else\the\rightheadtoks}%
    \ifsyntax@\setboxz@h{\def\\{\unskip\space\ignorespaces}%
        \headlinefont@##1}\fi}%
  \FN@\next@}
\def\NoRunningHeads{\global\runheads@false\global\let\headmark\eat@}

\newif\iffirstpage@     \firstpage@true
\newif\ifrunheads@      \runheads@true

\newdimen\fullhsize \fullhsize=\hsize
\newdimen\fullvsize \fullvsize=\vsize
\def\fullline{\hbox to\fullhsize}

\def\pagenumbers{\gdef\folio{\folio@}}

\let\norunningheads\NoRunningHeads
\def\userunningheads{\global\runheads@true}
\norunningheads

\headline={\def\chapter#1{\chapterno@. }%
  \def\\{\unskip\space\ignorespaces}\ifrunheads@\headlinefont@
    \ifodd\pageno\rightheadline \else\leftheadline\fi
   \else\hfil\fi\ifNoRunHeadline\global\NoRunHeadlinefalse\fi}
\let\folio@\folio
\def\foliofont@{\foliofont}
\def\foliofont{\eightrm}
\def\headlinefont@{\headlinefont}
\def\headlinefont{\eightpoint\smc}
\def\leftheadline{\rlap{\folio}\hfill
   \ifNoRunHeadline\else\iftrue\topmark\fi\fi \hfill}
\def\rightheadline{\hfill\ifNoRunHeadline
   \else \expandafter\fi
  \hfill \llap{\folio}}
\footline={{\eightpoint\bottremark}%
   \ifrunheads@\else\hfil{\let\foliofont\tenrm\folio}\fi\hfil}
\def\bottremark{}
 
\newif\ifNoRunHeadline      
\def\norunninghead{\global\NoRunHeadlinetrue}
\norunninghead

\output={\output@}
%
\newif\ifoffset\offsetfalse
\output={\output@}
\def\output@{%
 \ifoffset 
  \ifodd\count0\advance\hoffset by0.5truecm
   \else\advance\hoffset by-0.5truecm\fi\fi
 \shipout\vbox{%
  \makeheadline \pagebody \makefootline }%
 \advancepageno \ifnum\outputpenalty>-\@MM\else\dosupereject\fi}

\def\indexoutput#1{%
  \ifoffset 
   \ifodd\count0\advance\hoffset by0.5truecm
    \else\advance\hoffset by-0.5truecm\fi\fi
  \shipout\vbox{\makeheadline
  \vbox to\fullvsize{\boxmaxdepth\maxdepth%
  \ifvoid\topins\else\unvbox\topins\fi%
  #1 %
  \ifvoid\footins\else 
    \vskip\skip\footins
    \footnoterule
    \unvbox\footins\fi
  \ifr@ggedbottom \kern-\dimen@ \vfil \fi}%
  \baselineskip2pc
  \makefootline}%
 \global\advance\pageno\@ne
 \ifnum\outputpenalty>-\@MM\else\dosupereject\fi}
 
 \newbox\partialpage \newdimen\halfsize \halfsize=0.5\fullhsize
 \advance\halfsize by-0.5em

 \def\begindoublecolumns{\output={\indexoutput{\unvbox255}}%
   \begingroup \def\line{\fullline}
   \output={\global\setbox\partialpage=\vbox{\unvbox255\bigskip}}\eject
   \output={\doublecolumnout}\hsize=\halfsize \vsize=2\fullvsize}
 \def\enddoublecolumns{\output={\balancecolumns}\eject
  \endgroup \pagegoal=\fullvsize%
  \output={\output@}}
\def\doublecolumnout{\splittopskip=\topskip \splitmaxdepth=\maxdepth
  \dimen@=\fullvsize \advance\dimen@ by-\ht\partialpage
  \setbox0=\vsplit255 to \dimen@ \setbox2=\vsplit255 to \dimen@
  \indexoutput{\pagesofar} \unvbox255 \penalty\outputpenalty}
\def\pagesofar{\unvbox\partialpage
  \wd0=\hsize \wd2=\hsize \hbox to\fullhsize{\box0\hfil\box2}}
\def\balancecolumns{\setbox0=\vbox{\unvbox255} \dimen@=\ht0
  \advance\dimen@ by\topskip \advance\dimen@ by-\baselineskip
  \divide\dimen@ by2 \splittopskip=\topskip
  {\vbadness=10000 \loop \global\setbox3=\copy0
    \global\setbox1=\vsplit3 to\dimen@
    \ifdim\ht3>\dimen@ \global\advance\dimen@ by1pt \repeat}
  \setbox0=\vbox to\dimen@{\unvbox1} \setbox2=\vbox to\dimen@{\unvbox3}
  \pagesofar}

\tenpoint
\catcode`\@=\active

\def\smallheadings{\let\chapheadfonts\tenpoint\let\headfonts\tenpoint}

\tenpoint
\catcode`\@=\active

\magnification1200

\TagsOnRight

\def\KratBN{1}
\def\KratBZ{2}
\def\KuNiAA{3}
\def\MillAA{4}
\def\SadoAA{5}
\def\OEIS{6}
\def\EyndAA{7}

\def\TA{1}
\def\TB{2}
\def\TC{3}
\def\TD{4}
\def\TE{5}

\def\AA{2.1}
\def\AB{2.2}
\def\AC{2.3}
\def\AD{3.1}
\def\AE{3.2}
\def\AF{3.3}
\def\AG{3.4}
\def\AH{3.5}

\def\BA{5.1}
\def\BB{5.2}
\def\BC{5.3}
\def\BD{5.4}
\def\BF{5.5}
\def\BFa{5.6}
\def\BFb{5.7}
\def\BFc{5.8}
\def\BG{5.9}
\def\BH{5.10}

\def\CA{5.11}
\def\CB{5.12}

\def\DA{6.1}
\def\DAa{6.2}
\def\DBa{6.3}
\def\DBb{6.4}
\def\DBc{6.5}
\def\DBd{6.6}
\def\DBe{6.7}
\def\DBf{6.8}
\def\DBg{6.9}
\def\DC{6.10}
\def\DD{6.11}
\def\DE{6.12}
\def\DEa{6.13}
\def\DF{6.14}
\def\DG{6.15}
\def\DH{6.16}
\def\DHa{6.17}
\def\DI{6.18}
\def\DJ{6.19}
\def\DK{6.20}
\def\DL{6.21}

\catcode`\@=11
\def\iddots{\mathinner{\mkern1mu\raise\p@\hbox{.}\mkern2mu
    \raise4\p@\hbox{.}\mkern2mu\raise7\p@\vbox{\kern7\p@\hbox{.}}\mkern1mu}}
\catcode`\@=13

\def\zw#1{\left\langle#1\right\rangle}

\topmatter 
\title Spiral determinants
\endtitle 
\author G. Bhatnagar and C.~Krattenthaler 
\endauthor 
\affil 
Fakult\"at f\"ur Mathematik, Universit\"at Wien,\\
Oskar-Morgenstern-Platz~1, A-1090 Vienna, Austria.\\
WWW: {\tt http://www.mat.univie.ac.at/\~{}kratt}
\endaffil
\address Fakult\"at f\"ur Mathematik, Universit\"at Wien,
Oskar-Morgenstern-Platz~1, A-1090 Vienna,\linebreak Austria.\newline
WWW: \tt http://www.mat.univie.ac.at/\~{}kratt
\endaddress

\thanks Research supported 
by the Austrian Science Foundation FWF, 
grant SFB F50 (Special Research Program
``Algorithmic and Enumerative Combinatorics")
\endthanks
\subjclass Primary 15A15; 
 Secondary 05A10 11C20
\endsubjclass
\keywords spiral determinants
\endkeywords
\abstract 
We evaluate determinants of ``spiral" matrices, which are matrices in
which entries are spiralling from the centre of the matrices towards
the outside, with prescribed increments from one entry to the next
depending on whether one moves right, up, left, or down along the
spiral.
\endabstract
\endtopmatter
\document

\subhead 1. Preamble \endsubhead
The evaluation of determinants is a rich topic, and fascinating for many
people, in particular so for the authors of this note. A large body of
such evaluations has been collected in \cite{\KratBN, \KratBZ}.
Although it may seem to some that everything which is known is to be found
there, this is definitely not the case. After all, the task of assembling
``everything which is known" is an impossible one. (However, more is
contained in these two surveys as one might commonly think since frequently
simple transformations or reindexing turns a seemingly ``unknown"
determinant into a known one.) The purpose of this note is to discuss 
some attractive determinant evaluations not contained
(in any form) in \cite{\KratBN, \KratBZ}.

\subhead 2. Spiral matrices\endsubhead
Rather than jumping directly to our results in Theorems~\TA--\TC\
below (and the generalisation of Theorem~\TC\ discussed in Section~6)
without any further motivation,
we believe that it is of interest to describe the path that 
led to them.

During a visit of the second author (CK) at the Institut Henri
Poincar\'e in Paris in February 2017, Alexander R. Miller told
CK that he had been looking at ``spiral" matrices of the form
$$\pmatrix 
q^{16}&q^{15}&q^{14}&q^{13}\\
q^{5}&q^{4}&q^{3}&q^{12}\\
q^{6}&q&q^{2}&q^{11}\\
q^{7}&q^{8}&q^{9}&q^{10}
\endpmatrix\quad 
\text{and}\quad 
\pmatrix 
q^{17}&q^{16}&q^{15}&q^{14}&q^{13}\\
q^{18}&q^{5}&q^{4}&q^{3}&q^{12}\\
q^{19}&q^{6}&q&q^{2}&q^{11}\\
q^{20}&q^{7}&q^{8}&q^{9}&q^{10}\\
q^{21}&q^{22}&q^{23}&q^{24}&q^{25}
\endpmatrix,\quad \text{etc.}
\tag\AA
$$
(The pattern should be clear: we start with $q$ in the centre, and
then starting from there form a spiral of entries
$q,q^2,q^3,\dots,q^{n^2}$, the matrices above being the cases for
$n=4$ and $n=5$.) He had observed that the determinants of these
matrices all factor nicely, and he had also an explanation why
this happened. He asked CK whether he had already seen such
determinant evaluations before.

The reaction of CK was that, while he had not seen these determinants 
before, he seems to recall determinants of spiral matrices to have 
appeared in a Monthly Problem several years ago. It did not take long to
retrace this problem (namely \cite{\SadoAA}, which, upon further
search, turned out to have appeared even earlier as a problem in
the Mathematics Magazine \cite{\EyndAA}). 
As it turned out, it concerned the spiral matrices
$$\pmatrix 
1&2&3&4\\
12&13&14&5\\
11&16&15&6\\
10&9&8&7
\endpmatrix\quad 
\text{and}\quad 
\pmatrix 
1&2&3&4&5\\
16&17&18&19&6\\
15&24&25&20&7\\
14&23&22&21&8\\
13&12&11&10&9
\endpmatrix,\quad \text{etc.,}
\tag\AB
$$
where $1,2,\dots,n^2$ is winding from the outside into the centre.
Again, as the solution \cite{\MillAA} showed, the determinants of
all these matrices are given by nice product formulae.

Miller then went on to ask himself what happens if, instead of
spiralling from outside to inside, we spiral from inside to outside.
In other words, what are the determinants of the matrices
$$\pmatrix 
{16}&{15}&{14}&{13}\\
{5}&{4}&{3}&{12}\\
{6}&1&{2}&{11}\\
{7}&{8}&{9}&{10}
\endpmatrix\quad 
\text{and}\quad 
\pmatrix 
{17}&{16}&{15}&{14}&{13}\\
{18}&{5}&{4}&{3}&{12}\\
{19}&{6}&1&{2}&{11}\\
{20}&{7}&{8}&{9}&{10}\\
{21}&{22}&{23}&{24}&{25}
\endpmatrix,\quad \text{etc.?}
\tag\AC
$$
The On-Line Encyclopedia of Integer Sequences \cite{\OEIS} told him that the
determinants of these matrices were given by sequence {\tt A079340},
and a conjectured formula could be found there. (The determinants of
the matrices in (\AB) form the sequence {\tt A023999} in \cite{\OEIS}.)
He reported this to CK and asked him if he would know how to evaluate
the determinants in (\AB) and~(\AC).

After a few moments' thought,
CK's reply was that, say for the case of odd~$n$ ($n$ being the number
of rows and columns of the matrix), one should subtract the
next-to-last row from the last; the result is that almost all entries
in the last row are the same, except for the one in the first column.
Now use the first column to eliminate all the other entries in the
last row and expand the determinant with respect to the (new) last
row. As a result, one is again left with the determinant of a spiral 
matrix. Doing this, and something similar in the case of even~$n$, 
should result into an inductive proof.

When CK thought this argument through more carefully, he discovered
that, while it is true that the above manipulations allow one to
reduce the evaluation of the determinants of the spiral matrices in
(\AB) and (\AC) to the evaluation of a determinant of a smaller
spiral matrix, this smaller spiral matrix does however not belong
to the same family anymore. Neither is the central entry still
$1$ respectively~$n^2$, nor are the increments between entries 
along the spiral anymore always~1. Thus, if an inductive argument
in this spirit should work at all, then more parameters would have to
be introduced into the game.

\subhead 3. More parameters\endsubhead
Indeed, to introduce more parameters is one of the important points
``preached" in \cite{\KratBN, \KratBZ}: one should always try to
add parameters to the determinants one considers, and hope that
even with these new parameters the determinants still evaluate nicely.
If this happens, it will (usually) be much easier to evaluate the
more general determinants since parameters allow a much bigger
flexibility in the arguments.

In our case, what should we do? First of all, instead of~1 (or~$n^2$), 
we should allow ourselves to start with an arbitrary number in the
centre --- say~$a$ ---, and we should allow more general increments
from one entry to the next along the spiral. Clearly, if one allows
arbitrary increments then one cannot expect anything (since one then
faces a completely generic matrix), but maybe we should allow 
an (additive) increment of~$x$, say, when moving to the right,
an increment of~$b$, say, when moving up,
an increment of~$y$, say, when moving to the left,
and an increment of~$c$, say, when moving down along the spiral.
We denote the corresponding $n\times n$ spiral
matrix by $M_n(a,b,c,x,y)$. When it is clear from the context what
the parameters are, we shall abbreviate this to $M_n$. For example,
we have

\vskip-4pt
{\eightpoint
$$M_4=\pmatrix 
a+4 b+2 c+4 x+5 y & a+4 b+2 c+4 x+4 y & a+4 b+2 c+4 x+3 y & a+4 b+2 c+4 x+2
   y \\
 a+b+x+2 y & a+b+x+y & a+b+x & a+3 b+2 c+4 x+2 y \\
 a+b+c+x+2 y & a & a+x & a+2 b+2 c+4 x+2 y \\
 a+b+2 c+x+2 y & a+b+2 c+2 x+2 y & a+b+2 c+3 x+2 y & a+b+2 c+4 x+2 y 
\endpmatrix,
$$}%
and 
$$
M_5=\left(\smallmatrix 
a+4 b+2 c+4 x+6 y & a+4 b+2 c+4 x+5 y & a+4 b+2 c+4 x+4 y & a+4 b+2 c+4 x+3
   y & a+4 b+2 c+4 x+2 y \\
 a+4 b+3 c+4 x+6 y & a+b+x+2 y & a+b+x+y & a+b+x & a+3 b+2 c+4 x+2 y \\
 a+4 b+4 c+4 x+6 y & a+b+c+x+2 y & a & a+x & a+2 b+2 c+4 x+2 y \\
 a+4 b+5 c+4 x+6 y & a+b+2 c+x+2 y & a+b+2 c+2 x+2 y & a+b+2 c+3 x+2 y &
   a+b+2 c+4 x+2 y \\
 a+4 b+6 c+4 x+6 y & a+4 b+6 c+5 x+6 y & a+4 b+6 c+6 x+6 y & a+4 b+6 c+7 x+6
   y & a+4 b+6 c+8 x+6 y 
\endsmallmatrix\right).
$$

Computer experiments suggest that
the determinants of these matrices factorise nicely.
In fact, it is not difficult to come up with a guess for the result,
and the inductive argument sketched above leads to an almost effortless proof.

\proclaim{Theorem \TA}
For all non-negative integers $n$, we have
$$\multline 
\det M_{2n}(a,b,c,x,y)=(-1)^{n+1}
\big(a x + n^2 b x + n (n - 1) c x + n^2 x^2 + a y + (n - 1)^2 b y \\
+ 
 n (n - 1) c y + n (n - 1) y^2 + n (2 n - 1) x y\big)
\prod _{i=1} ^{2n-2}\big(i(b+c)+(i+1)(x+y)\big)
\endmultline$$
and
$$\multline 
\det M_{2n+1}(a,b,c,x,y)=(-1)^{n}
\big(a x + n^2 b x + n (n - 1) c x + n^2 x^2 + a y + n^2 b y \\
+ 
 n (n + 1) c y + n (n + 1) y^2 + n (2 n + 1) x y\big)
\prod _{i=1} ^{2n-1}\big(i(b+c)+(i+1)(x+y)\big).
\endmultline$$
\endproclaim

\demo{Proof}
One proceeds by induction on~$n$.
Let us consider the matrix $M_{2n+1}(a,b,c,x,y)$. It has the form
$$
M_{2n+1}=
\pmatrix 
\hdotsfor5\\
\vdots& & & &\vdots\\
E_1&\dots&a&\innerhdotsfor2\after\quad &\\
\hdotsfor5\\
E_2&E_3&\innerhdotsfor2\after\quad &E_4\\
E_5&E_6&\innerhdotsfor2\after\quad &E_7
\endpmatrix,
\tag\AD
$$
where
$$\align 
E_1&=a+n^2b+n^2c+n^2x+n(n+1)y,\\
E_2&=a+n^2b+(n^2+n-1)c+n^2x+n(n+1)y,\\
E_3&=a+(n-1)^2b+n(n-1)c+(n-1)^2x+n(n-1)y,\\
E_4&=a+(n-1)^2b+n(n-1)c+n^2x+n(n-1)y,\\
E_5&=a+n^2b+n(n+1)c+n^2x+n(n+1)y,\\
E_6&=a+n^2b+n(n+1)c+(n^2+1)x+n(n+1)y,\\
E_7&=a+n^2b+n(n+1)c+n(n+2)x+n(n+1)y.
\endalign$$

Now we subtract the next-to-last row from the last row.
This operation does not change the determinant. The matrix
we obtain thereby is
$$
\pmatrix 
\hdotsfor5\\
\vdots& & & &\vdots\\
E_1&\dots&a&\innerhdotsfor2\after\quad &\\
\hdotsfor5\\
E_2&E_3&\innerhdotsfor2\after\quad &E_4\\
c&D_1&\innerhdotsfor2\after\quad &D_1
\endpmatrix,$$
where
$$
D_1=(2n-1)b+2nc+2nx+2ny.
$$
We eliminate the $D_1$'s by subtracting $D_1/c$ times the first column
from the other columns. Again, this does not change the determinant.
The resulting matrix has the form
$$
\pmatrix 
\hdotsfor5\\
\vdots& & & &\vdots\\
E_1&\dots&A_1&\innerhdotsfor2\after\quad &\\
\hdotsfor5\\
*&*&\innerhdotsfor2\after\quad &*\\
c&0&\innerhdotsfor2\after\quad &0
\endpmatrix,
\tag\AE$$
where
$$\align 
A_1&=a-\frac {1} {c}D_1E_1\\
&=a-\frac {1} {c}\big((2n-1)b+2nc+2nx+2ny\big)
\big(a+n^2b+n^2c+n^2x+n(n+1)y\big).
\tag\AF
\endalign$$
The point here is that, apart from the first column and the last row,
we have obtained a spiral matrix, of dimensions $(2n)\times(2n)$,
with central entry~$A_1$, with right increment~$x$, up increment~$B_1$,
left increment~$y$, and down increment~$C_1$, where
$$\align 
B_1&=b+D_1=2n(b+c+x+y),
\tag\AG\\
C_1&=c-D_1=-(2n-1)(b+c)-2n(x+y).
\tag\AH
\endalign$$
Hence, by taking the determinants of (\AD) and (\AE), we infer
the relation
$$
\det M_{2n+1}(a,b,c,x,y)
=
c\cdot\det M_{2n}(A_1,B_1,C_1,x,y),
$$
where $A_1,B_1,C_1$ are given in (\AF)--(\AH).
If we assume by induction that we already know the formula for $M_{2n}$,
then by the above relation we obtain the claimed formula for $M_{2n+1}$.

For the proof of the claim on the determinant of $M_{2n}(a,b,c,x,y)$,
one proceeds similarly. Here, one subtracts the second row from the
first row, and subsequently $D_2/b$ times the last column from the
other columns, where
$$
D_2=(2n-1)b+(2n-2)c+(2n-1)x+(2n-1)y.
$$
Here, the relation which results by taking determinants is
$$
\det M_{2n}(a,b,c,x,y)
=
-b\cdot\det M_{2n-1}(A_2,B_2,C_2,x,y),
$$
where
$$\align 
A_2&=a-\frac {1} {b}\big((2n-1)b+(2n-2)c+(2n-1)x+(2n-1)y\big)\\
&\kern2cm\times
\big(a+n(n-1)b+n(n-1)c+n^2x+n(n-1)y\big),\\
B_2&=-(2n-2)(b+c)-(2n-1)(x+y),\\
C_2&=(2n-1)(b+c+x+y).
\endalign$$
Again, if we assume by induction that we already know the formula
for $M_{2n-1}$, then the claimed formula for $M_{2n}$ results
straightforwardly from the above relation.\quad \quad \qed
\enddemo

\remark{Remark}
It goes without saying that the evaluation of the determinants of the
matrices in (\AB) results by specialising Theorem~\TA\ to $a=n^2$,
$x=y=b=c=-1$, while the evaluation of the determinants of the matrices
in (\AC) is given by the special case $a=x=y=b=c=1$ of Theorem~\TA.
In particular, this provides a proof of the conjectured formula in
\cite{\OEIS, {\tt A079340}}. The proof \cite{\MillAA} of the formula
for the determinants in (\AB) actually yielded more generally the
evaluation of the determinants of the matrices $M_n(a,-1,-1,-1,-1)$.
\endremark

\subhead 4. A ``$q$-analogue"\endsubhead
The ``$q$-analogue" of the above theorem in which one replaces each
entry $X$ of $M_n(a,b,c,x,y)$ by $q^X$ is actually simpler to
derive. Let $Q_n(a,b,c,x,y)$ denote the matrix which results from
this replacement. Clearly, the matrices in (\AA) are 
the matrices $Q_n(1,1,1,1,1)$.

\proclaim{Theorem \TB}
For all non-negative integers $n$, we have
$$\multline 
\det Q_{2n}(a,b,c,x,y)=(-1)^{n}
q^{2n a + \frac {1} {3}n (2 n^2 + 1) (b + x) + 
     \frac {2} {3} (n - 1) n (n + 1) (c + y)}\\
\times
\prod _{i=0} ^{2n-2}\big(1-q^{i(b+c)+(i+1)(x+y)}\big)
\endmultline$$
and
$$\multline
\det Q_{2n+1}(a,b,c,x,y)=(-1)^{n}
  q^{(2n+1) a + \frac {1} {3}n (n + 1) (2 n + 1) (b + c + x + y)}\\
\times
\prod _{i=0} ^{2n-1}\big(1-q^{i(b+c)+(i+1)(x+y)}\big).
\endmultline$$
\endproclaim

\demo{Proof}
One proceeds analogously to before.
When one does the subtraction of the rows (where one
multiplies by a suitable power of $q$), then all entries but one
in the last/first row already become zero. So, no column operations
are needed. Everything else is just the same.\quad \quad \qed
\enddemo
\remark{Remark}
The evaluation of the determinants of the
matrices in (\AA) results by specialising Theorem~\TB\ to 
$a=x=y=b=c=1$.
\endremark

\subhead 5. A real $q$-analogue\endsubhead
Alexander Miller pointed out that the ``$q$-analogue"
given in Theorem~\TB\ is not a {\it real\/} $q$-analogue: no matter how
one specialises Theorem~\TB\ or which limit one takes, it will never be
possible to derive Theorem~\TA.

Without doubt, this is a valid point. 
How can one obtain a ``real" $q$-analogue? 
The standard way (in combinatorics) 
to pass from the ``ordinary
$q=1$-world" to the ``$q$-world" is by replacing entries~$X$ by
$q^X-1$. If one tries this with the matrix $M_n(a,b,c,x,y)$ in
Theorem~\TA, then the determinant of the resulting matrix 
turns out to be a huge polynomial
in $q^a,q^b,q^c,q^x,q^y$ which does not factor at all. 

A second way to arrive at $q$-analogues is to replace entries~$X$
by the symmetric $q^X-q^{-X}$. (This is very common in the theory
of quantum algebras.) Again the determinants of the resulting matrices
 are  huge expressions that do not factor at all. 
 However, on setting $x=y$, a miracle happens, and we obtain 
very nice factorisations. (For the record: nothing of this kind happens if one
sets $b=c$, or for any other similar specialisation.) Again, it is
not difficult to guess the result. See Theorem~\TC\ below. 

The recursive approach which proved Theorems~\TA\ and~\TB\ does not work
anymore. However, by calculating a few small cases ---
this time by hand rather than by a computer --- 
an elimination technique was found which led to a proof.



In order to formulate the result and its proof in a convenient and
succinct way, we ``get rid" of~$q$: instead of writing $q^X-q^{-X}$,
we write $X-\frac {1} {X}$. To be precise,
let $Z_n(a,b,c,x,y)$ be the $n\times n$ spiral matrix, where, instead
of adding increments of $x,b,y,c$ at each ``step" when
moving to the right, up, left, down in the spiral, we {\it multiply}
by $x,b,y,c$, and in the end, if the result is $\alpha$, we put
$\alpha-\alpha^{-1}$ as the corresponding entry.
For the sake of simplicity
of notation, we write $[\alpha]:=\alpha-\alpha^{-1}$ and $\zw{\alpha}:=\alpha+\alpha^{-1}$.
For example, the $4\times4$ matrix $Z_4(a,b,c,x,y)$ is
$$Z_4=\pmatrix 
[ab^4c^2x^4y^5 ]&[ ab^4c^2x^4y^4 ]&[ ab^4c^2x^4y^3 ]&[ ab^4c^2x^4y^2 ]\\
 [abxy^2 ]&[ abxy ]&[ abx ]&[ ab^3c^2x^4y^2] \\
 [abcxy^2 ]&[ a ]&[ ax ]&[ ab^2c^2x^4y^2] \\
 [abc^2xy^2 ]&[ abc^2x^2y^2 ]&[ abc^2x^3y^2 ]&[ abc^2x^4y^2 ]
\endpmatrix.
$$
In order to see that this is a $q$-analogue of $M_n(a,b,c,x,y)$, one has
to replace each variable $\alpha$ by $q^\alpha$, divide each entry by
$q-q^{-1}$, and then let $q\to1$.

The announced determinant evaluation is the following.

\proclaim{Theorem \TC}
For all non-negative integers $n$, we have
$$\multline 
\det Z_{2n}(a,b,c,x,x)=(-1)^{n+1}
\left[a^2b^{2n^2-2n+1}c^{2n^2-2n}x^{2n(2n-1)}\right]
\prod _{k=0} ^{2n-2}\left[(bc)^{k/2}x^{k+1}\right]\\
\times
\prod _{k=1} ^{n-1}\zw{ab^{k(k+1)}c^{k^2}x^{k(2k+1)}}
\zw{ab^{(2k^2-2k+1)/2}c^{(2k^2-1)/2}x^{k(2k-1)}}
\endmultline
\tag\BA$$
and
$$\multline 
\det Z_{2n+1}(a,b,c,x,x)=(-1)^{n}
\left[a^2(bc)^{2n^2}x^{2n(2n+1)}\right]
\prod _{k=0} ^{2n-1}\left[(bc)^{k/2}x^{k+1}\right]\\
\times
\prod _{k=1} ^{n-1}\zw{ab^{k(k+1)}c^{k^2}x^{k(2k+1)}}
\prod _{k=1} ^{n}\zw{ab^{(2k^2-2k+1)/2}c^{(2k^2-1)/2}x^{k(2k-1)}}.
\endmultline
\tag\BB$$
\endproclaim

\demo{Proof}
Our proof will be based on the two relations
$$\align 
\left[a\right] \langle x\rangle
&= \left[ax\right]+\left[a/x\right], 
\tag\BC\\
\left[a\right]\left[x\right]&= \left\langle ax\right\rangle
-\left\langle a/x\right\rangle ,
\tag\BD
\endalign$$
valid for all non-zero $a$ and $x$.
From (\BC), we obtain
$$
\left[Ax^2\right]-\zw{x}[Ax]+[A]=0,
\tag\BF
$$
for all non-zero $A$ and $x$, which is the key identity in the
following arguments.

\medskip
We start with the even case, that is, with the proof of (\BA).
In the following, we shall abbreviate $Z_{2n}(a,b,c,x,x)$ by
$Z_{2n}$. Furthermore, we denote the matrix entries by $z_{i,j}$,
so that $Z_{2n}=(z_{i,j})_{1\le i,j\le 2n}$.

Let $C_j$ denote the $j$-th column of $Z_{2n}$. We replace $C_j$
by
$$
C_j-\zw{x}C_{j-1}+C_{j-2},\quad \text{for }j=3,4,\dots,2n.
\tag\BFa
$$
Let $Z_{2n}'$ be the matrix which is obtained in this way.
Clearly, these column operations do not change the determinant,
so that
$$
\det Z_{2n}=\det Z_{2n}'.
\tag\BFb
$$
By (\BF), the effect of these operations is that, starting from the
third column, all entries of $Z_{2n}'$ in the
upper wedge bordered by the diagonal $z_{i,i}'$, $i=1,2,\dots,n$
and the antidiagonal $z_{i,2n-i+2}'$, $i=2,3,\dots,n$ vanish, and,
similarly, that all entries of $Z_{2n}'$ in the
lower wedge bordered by the diagonal $z_{i,i+1}'$, $i=n+1,n+2,\dots,2n-1$
and the antidiagonal $z_{i,2n-i+2}'$, $i=n+1,n+2,\dots,2n$ vanish.
More precisely, the matrix $Z_{2n}'$ has the form
$$
\pmatrix\format\c&\c&\c&\c&\c&\c&\c&\c&\c&\c&\c&\c\\
z_{1,1} & z_{1,2} & 0 & 0 & \innerhdotsfor5\after\quad & 0& 0 & 0\\
* & * & 0 & 0 & \innerhdotsfor5\after\quad & 0 & 0 & z_{2,2n}' \\ 
* & * & z_{3,3}' & 0 & \innerhdotsfor5\after\quad & 0 & z_{3,2n-1}' & * \\ 
* & * & * & z_{4,4}' & 0 & \innerhdotsfor3\after\quad & 0 & z_{4,2n-2}' & * & * \\ 
\vdots & \vdots & & & \ddots & & & & \iddots & &  &\vdots \\
\vdots & \vdots & & & & z_{n,n}' & 0 & z_{n,n+2}' & & &  &\vdots \\
\vdots & \vdots & & & & & z_{n+1,n+1}'\ & z_{n+1,n+2}' & & &  &\vdots \\
\vdots & \vdots & & & & \iddots & & & \ddots& & & \vdots\\
\vdots & \vdots & & & \iddots & & & & & \ddots& & \vdots\\
* & * & * & z_{2n-2,4}' & 0 & \innerhdotsfor4\after\quad & 0& z_{2n-2,2n-1}'& *\\
* & * & z_{2n-1,3}' & 0 & \innerhdotsfor5\after\quad & 0& 0 & z_{2n-1,2n}'\\
z_{2n,1}\ & z_{2n,2} & 0 & 0 & \innerhdotsfor5\after\quad & 0& 0 & 0
\endpmatrix .
\tag\BFc
$$

We may now do a Laplace expansion of the determinant of this matrix 
with respect to the first and the last row. The result is that
$$\multline
\det Z_{2n}'=\det
\pmatrix z_{1,1}&z_{1,2}\\
z_{2n,1}&z_{2n,2} \endpmatrix\\
\times
\det\pmatrix\format\c&\c&\c&\c&\c&\c&\c&\c&\c&\c&\c&\c\\
0 & 0 & \innerhdotsfor5\after\quad & 0 & 0 & z_{2,2n}' \\ 
z_{3,3}' & 0 & \innerhdotsfor5\after\quad & 0 & z_{3,2n-1}' & * \\ 
* & z_{4,4}' & 0 & \innerhdotsfor3\after\quad & 0 & z_{4,2n-2}' & * & * \\ 
\vdots & & \ddots & & & & \iddots & &  &\vdots \\
\vdots & & & z_{n,n}' & 0 & z_{n,n+2}' & & &  &\vdots \\
\vdots & & & & z_{n+1,n+1}'\ & z_{n+1,n+2}' & & &  &\vdots \\
\vdots & & & \iddots & & & \ddots& & & \vdots\\
\vdots & & \iddots & & & & & \ddots& & \vdots\\
* & z_{2n-2,4}' & 0 & \innerhdotsfor4\after\quad & 0& z_{2n-2,2n-1}'& *\\
z_{2n-1,3}' & 0 & \innerhdotsfor5\after\quad & 0& 0 & z_{2n-1,2n}'
\endpmatrix .
\endmultline
$$
By expanding the last determinant alternatingly along the first row,
along the last row, along the second row, along the next-to-last row,
etc., it is straightforward to see that --- up to a sign that is
easily determined --- 
it equals the product of the entries along the antidiagonal
$z_{i,2n-i+2}'$, $i=2,3,\dots,2n-1$. Thus, we obtain
$$
\det Z_{2n}'=(-1)^{n-1}(z_{1,1}z_{2n,2}-z_{1,2}z_{2n,1})
\prod _{i=2} ^{2n-1}z_{i,2n-i+2}'.
\tag\BG
$$
Hence, the only remaining task is to compute the first term in
parentheses on the right-hand side and the entries $z_{i,2n-i+2}'$.

We begin with the first term on the right-hand side of (\BG):
$$\align 
z_{1,1}z_{2n,2}-z_{1,2}z_{2n,1}
&=
\left[ a b^{n^2} c^{n(n-1)} x^{n(2n-1)+(2n-1)}\right]
\left[ a b^{(n-1)^2} c^{n(n-1)} x^{n(2n-1)-(2n-2)}\right]\\
&\kern.5cm
-
\left[ a b^{n^2} c^{n(n-1)} x^{n(2n-1)+(2n-2)}\right]
\left[ a b^{(n-1)^2} c^{n(n-1)} x^{n(2n-1)-(2n-1)}\right]\\
& = \left\langle a^2 b^{2n^2-2n+1} c^{2n^2-2n}
x^{2n(2n-1)+1}\right\rangle \\  
& \kern1cm -
\left\langle a^2 b^{2n^2-2n+1} c^{2n^2-2n}
x^{2n(2n-1)-1}\right\rangle \\ 
& = \left[ a^2 b^{2n^2-2n+1} c^{2n^2-2n} x^{2n(2n-1)}\right] \left[
x\right]  .
\tag\BH
\endalign$$
Here, we first used (\BD) to rewrite the products of the form
$[A][B]$, and subsequently we used (\BD) in the ``reverse" direction
in order to obtain the last line.

To complete the proof in the current case, we need a formula for the
entries $z_{i,2n-i+2}'$ in (\BG).  
To do that, we need formulae for the entries 
$z_{i,2n-i},z_{i,2n-i+1},z_{i,2n-i+2}$ of the original matrix $Z_{2n}$.
Writing $k=n-i+1$ for convenience, we find the formulae
$$\alignat2
z_{i,2n-i} &=\left[ a b^{k^2} c^{k(k-1)} x^{k(2k-1)+1}\right], \quad 
&&\text{for } i=1,2, \dots, n; 
\\
z_{2n-i+1, i-1} &=\left[ a b^{k^2} c^{k(k+1)-1}
x^{k(2k+1)}\right],\quad  &&\text{for } i=2,3, \dots, n;\\
z_{i,2n-i+1} &=\left[ a b^{k^2} c^{k(k-1)} x^{k(2k-1)}\right],\quad 
&&\text{for } i=1,2, \dots, n; 
\\
z_{2n-i+1, i} &=\left[ a b^{(k-1)^2} c^{k(k-1)}
x^{(k-1)(2k-1)}\right],
\quad  &&\text{for } i=1,2, \dots, n;\\
z_{i,2n-i+2} &=\left[ a b^{k^2+2k} c^{k(k+1)} x^{k(2k+1)}\right]
,\quad  &&\text{for } i=2,3, \dots, n; 
\\
z_{2n-i+1, i+1} &=\left[ a b^{(k-1)^2} c^{k(k-1)}
x^{(k-1)(2k-1)+1}\right],\quad && \text{for } i=1,2, \dots, n .
\endalignat
$$
Thus, by using (\BC) several times, for $i= 2, 3, \dots, n$ we have
$$\align
z_{i,2n-i+2}' &= z_{i,2n-i+2} - \langle x\rangle z_{i,2n-i+1} 
+ z_{i,2n-i} \\ 
&=  
\left[ a b^{k^2+2k} c^{k(k+1)} x^{(k+1)(2k+1)}\right] -
\langle x\rangle\left[ a b^{k^2} c^{k(k-1)} x^{k(2k-1)}\right]
 \\
& \kern1cm + \left[ a b^{k^2} c^{k(k-1)} x^{k(2k-1)+1}\right] \\
& = \left[ a b^{k^2+2k} c^{k(k+1)} x^{(k+1)(2k+1)}\right] -
\left[ a b^{k^2} c^{k(k-1)} x^{k(2k-1)-1}\right]
 \\
& = \left[ (bc)^{k}  x^{2k+1}\right]
\left\langle a b^{k(k+1)} c^{k^2} x^{k(2k+1)}\right\rangle
\tag\CA
\endalign
$$
and
$$\align
z_{2n-i+1,i+1}' &= z_{2n-i+1,i+1} - \langle
x\rangle z_{2n-i+1,i}  + z_{2n-i+1, i-1} \\ 
&=  
\left[ a b^{(k-1)^2} c^{k(k-1)} x^{(k-1)(2k-1)+1}\right] -
\langle x\rangle\left[ a b^{(k-1)^2} c^{k(k-1)} x^{(k-1)(2k-1)}\right]
 \\
& \kern1cm + \left[ a b^{k^2} c^{k(k+1)-1} x^{k(2k+1)}\right] \\
& = -\left[ a b^{(k-1)^2} c^{k(k-1)} x^{(k-1)(2k-1)}\right]
+\left[ a b^{k^2} c^{k(k+1)-1} x^{k(2k+1)}\right]
 \\
& = \left[ (bc)^{(2k-1)/2}  x^{2k}\right]
\left\langle a b^{(2k^2-2k+1)/2} c^{(2k^2-1)/2} x^{k(2k-1)}\right\rangle.
\tag\CB
\endalign
$$
Substitution of (\BH), (\CA), and (\CB) in (\BG) in combination with
(\BFb) leads to (\BA).

\medskip
The odd case, that is, our proof of (\BB), is similar.
The only difference is that, here, we replace the $j$-th column
of $Z_{2n+1}(a,b,c,x,x)$, $C_j$ say, by
$$C_j - \langle x \rangle C_{j+1} +C_{j+2},\quad \text{for }
j=1, 2, \dots, 2n-1.
$$ 
(This should be compared with (\BFa).)
Everything else is completely analogous. We leave the details to the
reader.

This completes the  proof of the theorem.\quad \quad \qed 
\enddemo

\subhead 6. Discussion \endsubhead
We may ask ourselves whether Theorem~\TC\ is the most general
result existing or whether further generalisations are possible.
In particular, Theorem~\TC\ is not a full generalisation of Theorem~\TA\
but only of its special case in which $x=y$. Clearly, our proof of
Theorem~\TC\ heavily depends on~$x$ being equal to~$y$: the column
operations that we perform there create zeroes in the upper {\it and\/} lower
wedges only in this special case. We must leave it as an open
problem whether there is a {\it full} ``$q$-analogue" of Theorem~\TA\
(in the sense discussed at the very beginning of this section).

\medskip
Is it possible to generalise Theorem~\TC\ in other directions?
An attentive reader will have noticed that
our proof only hinged upon the two relations (\BC) and (\BD)
satisfied by the brackets $[\,.\,]$
(also involving the functions $\zw{\,.\,}$). An even closer perusal of
the proof shows that actually relation (\BC) is sufficient for a
factorisation of the determinant to happen. Namely, it implies
(\BF), which in turn makes the column operations (\BFa) create
the wedges of zeroes in the matrix $Z_{2n}'$ (cf\. (\BFc)).
Relation (\BD) is only used in the simplification of the
$2\times2$-determinant in (\BH). It then becomes clear that 
we could have even chosen {\it arbitrary} increments while moving
up or down one step along the spiral, and we would still obtain
a closed form product formula for the determinant of the corresponding
matrix. Moreover, entries which are in the left and the right
wedges do not influence the result at all.

The final question which we address is whether there are
more general functions (instead of the brackets $[\,.\,]$ and angle
brackets $\zw{\,.\,}$)?
So, are our functions $[\,.\,]$ and $\zw{\,.\,}$ the most general functions
possible which satisfy (\BC)?

Our final proposition --- which is based on an auxiliary result that
we state and prove separately as Lemma~\TE\ below --- 
says that, while there are in fact slightly more
general functions satisfying (\BC), they do not allow for an ``essential"
generalisation of Theorem~\TC\ along this line. See the more
detailed discussion of the significance of Proposition~\TD\ 
in the remark after the statement of the
proposition.

\proclaim{Proposition \TD}
Let $f$ and $g$ be continuous 
functions from $\Bbb R_{>0}$ to $\Bbb R$ satisfying the relation
$$
f(a)g(x)=f(ax)+f(a/x)
\tag\DA
$$
for all positive real numbers $a$ and $x$. Then there are three possibilities:

\roster 
\item $f(x)=0$,\quad  for $x\in\Bbb R_{>0}$;
\item 
$
f(x)=c_1+c_2\log x\quad \text{and}\quad
g(x)=2,\quad \text{for }x\in\Bbb R_{>0},
$\newline
for suitable real numbers $c_1$ and $c_2$;
\item  
$
f(x)=c_1x^\alpha+c_2x^{-\alpha}\quad \text{and}\quad 
g(x)=x^\alpha+x^{-\alpha},\quad \text{for }x\in\Bbb R_{>0},
$\newline
for a suitable real or (purely) imaginary number
$\alpha\ne0$, and for suitable real numbers $c_1$ and $c_2$.
\endroster
\endproclaim

\remark{Remark}
We recall that this proposition in combination with the proof of
Theorem~\TC\ says that, if we form the spiral matrix which arises
from $Z_n(a,b,c,x,x)$ by replacing the bracket $[\,.\,]$ by one
of the functions~$f(\,.\,)$ given in the proposition, then its
determinant will factor in closed product form.

Clearly, the first possibility is without interest, since it
produces the zero matrix. As is not difficult to see, 
the second possibility leads --- after reparametrisation --- 
to the earlier spiral matrices
$M_n(a,b,c,x,y)$ with $x=y$. By Theorem~\TA, we know that their determinants
factor nicely, but this happens even for the spiral
matrices $M_n(a,b,c,x,y)$ with {\it arbitrary} $x$ and $y$. 

The third possibility does indeed provide functions that are more
general than the brackets $[\,.\,]$, which define spiral
matrices whose determinants still factor nicely. For the sake of
brevity, we omit to work out the explicit statement of the 
corresponding result, which we leave to the reader.
Suffice it to say that it still does not allow for a full $q$-analogue
of Theorem~\TA, in the sense that a specialisation or limit would be
equivalent to Theorem~\TA, with $x$ and $y$ {\it different}. 
\endremark

\demo{Proof of Proposition \TD}
Clearly, $f\equiv0$ solves (\DA), which corresponds to possibility~(1).
Therefore, from now on, we assume that $f$ is not identically
zero. 

\medskip
We start with the observation that under the replacement $x\to 1/x$
nothing changes on the right-hand side of (\DA), and therefore our
assumption on~$f$ implies
$$
g(y)=g(1/y),\quad \text{for }y\in\Bbb R_{>0}.
\tag\DAa
$$

Next we use (\DA) with $x=y$, with $x=y^2$, and with $x=y^3$ 
to produce the system of equations
$$\align 
f(ya)g(y)&=f(y^2a)+f(a),
\tag\DBa\\
f(y^2a)g(y)&=f(y^3a)+f(ya),
\tag\DBb\\
f(y^3a)g(y)&=f(y^4a)+f(y^2a),
\tag\DBc\\
f(y^4a)g(y)&=f(y^5a)+f(y^3a),
\tag\DBd\\
f(y^5a)g(y)&=f(y^6a)+f(y^4a),
\tag\DBe\\
f(y^2a)g(y^2)&=f(y^4a)+f(a).
\tag\DBf\\
f(y^3a)g(y^3)&=f(y^6a)+f(a).
\tag\DBg
\endalign$$
By adding both sides of (\DBa) and (\DBc), we obtain
$$
g(y)\big(f(ya)+f(y^3a)\big)=f(a)+f(y^4a)+2f(y^2a).
$$
If we now use (\DBb) on the left-hand side and (\DBf) on the right-hand
side, then we arrive at the equation
$$
g^2(y)f(y^2a)=\big(g(y^2)+2\big)f(y^2a).
$$
By our assumption on~$f$, this implies the equation
$$
g(y^2)=g^2(y)-2,\quad \text{for }y\in\Bbb R_{>0}.
\tag\DC
$$

Similarly, if we add both sides of (\DBa) and (\DBe), respectively
both sides of (\DBb) and (\DBd), then we obtain
$$
\align
g(y)\big(f(ya)+f(y^5a)\big)&=f(a)+f(y^6a)+f(y^2a)+f(y^4a),\\
g(y)\big(f(y^2a)+f(y^4a)\big)&=f(ya)+f(y^5a)+2f(y^3a).
\endalign
$$
We use the second equation to rewrite the left-hand side of the first
equation, and (\DBg) to rewrite the right-hand side. The result is
$$
g(y)\big(g(y)\big(f(y^2a)+f(y^4a)\big)-2f(y^3a)\big)
=g(y^3)f(y^3a)+f(y^2a)+f(y^4a).
$$
Finally, we use (\DBc) on both sides and arrive at
$$
g(y)\big(g^2(y)f(y^3a)-2f(y^3a)\big)
=g(y^3)f(y^3a)+g(y)f(y^3a).
$$
By our assumption on~$f$, this implies the equation
$$
g(y^3)=g^3(y)-3g(y),\quad \text{for }y\in\Bbb R_{>0}.
\tag\DD
$$

By Lemma~\TE, we see that (\DAa), (\DC), and (\DD) 
together imply already the claim
on the form of the function~$g$, namely that
$g(x)=x^\alpha+x^{-\alpha}$ for a suitable real or (purely) imaginary number~$\alpha$.
(This includes possibility~(2), in which $g(x)\equiv2$, 
by choosing $\alpha=0$.)

It remains to solve the functional equation
$$
\left(x^\alpha+x^{-\alpha}\right)f(a)
=f(ax)+f(a/x),\quad \text{for }a,x\in \Bbb R_{>0}.
\tag\DE
$$

\medskip
We treat the case where $\alpha\ne0$ first.
It is trivial to check that $f(a)=a^\alpha$ and $f(a)=a^{-\alpha}$
are solutions of (\DE). 

Let $f$ be an arbitrary solution of (\DE). We determine $c_1$ and $c_2$ as
solutions of the system of linear equations
$$\align 
f(1)&=c_1+c_2,\\
f(2)&=c_12^\alpha+c_22^{-\alpha}.
\endalign$$
It should be noted that this system can indeed always be
(uniquely) solved since $\alpha\ne0$. Using (\DE) with $x=2$, 
it may then inductively be inferred that
$$
f(2^n)=c_1(2^n)^\alpha+c_2(2^n)^{-\alpha},\quad \text{for }n\in\Bbb Z.
\tag\DEa
$$
Furthermore, given $y$ and $z$ with 
$$
f(y)=c_1y^\alpha+c_2y^{-\alpha}\quad \text{and}\quad 
f(z)=c_1z^\alpha+c_2z^{-\alpha},
$$
Equation (\DE) with $a=\sqrt{yz}$ and $x=\sqrt{y/z}$ implies that
$$
f(\sqrt{yz})=c_1\big(\sqrt{yz}\big)^\alpha+c_2\big(\sqrt{yz}\big)^{-\alpha}.
$$
Iteration of this argument starting from the values in (\DEa) yields
a dense set of values $x$ for which $f(x)=c_1x^\alpha+c_2x^{-\alpha}$.
By continuity of $f$, the same equation must hold for {\it all\/}
$x\in\Bbb R_{>0}$. We thus arrived at possibility~(3).

\medskip
If $\alpha=0$, then (\DE) simplifies to
$$
2f(a)
=f(ax)+f(a/x),\quad \text{for }a,x\in \Bbb R_{>0}.
$$
Here, it is straightforward to see that $f(a)=1$ and $f(a)=\log a$ are
solutions of this equation. By an argument that is completely
analogous to the one before that led to possibility~(3), we 
conclude that, here, we arrive at possibility~(2).

\medskip
This completes the proof of the proposition.\quad \quad \qed
\enddemo

\proclaim{Lemma \TE}
Let $g:\Bbb R^{>0}\to\Bbb R$ be a continuous function satisfying the relations
$$\align
g(x)=g(1/x),\quad \text{for }x\in\Bbb R_{>0},
\tag\DF\\
g(x^2)=g^2(x)-2,\quad \text{for }x\in\Bbb R_{>0},
\tag\DG\\
g(x^3)=g^3(x)-3g(x),\quad \text{for }x\in\Bbb R_{>0}.
\tag\DH
\endalign$$
Then 
$$
g(x)=x^\alpha+x^{-\alpha}
=2\cos\big(i\alpha\log x\big),
\quad \text{for }x\in\Bbb R_{>0},
$$
for a suitable real or (purely) imaginary number $\alpha$.
\endproclaim

\remark{Remark}
A corollary is that $g$ satisfies
$$
g(a)g(x)=g(ax)+g(a/x),\quad \text{for }a,x\in\Bbb R_{>0}.
\tag\DHa
$$
In fact, one can show that there holds a(n almost) reverse statement:
if $g$ satisfies (\DHa), then either $g$ vanishes identically or $g$ satisfies
(\DF)--(\DH).
\endremark

\demo{Proof of Lemma \TE}
From (\DG) and (\DH) we see that $g(1)=2$.
Because of (\DF), we may restrict our attention to the values $g(x)$
with $x\ge1$, which we shall do from now on.

Our strategy is to find an $\alpha$ and a set $\{ x_{n,m} : n, m = 0, 1,
2, \dots  \}$  that is dense in $\Bbb R_{>0}$ such that
$$
g(x_{n,m})=x_{n,m}^\alpha+x_{n,m}^{-\alpha},\quad 
\text{for }n,m=0,1,2,\dots.
\tag\DI
$$

We distinguish between two cases: 
\roster 
\item there exists $x_0>1$ such that $g(x_0)\ge2$;
\item there exists $x_0>1$ such that $g(x_0)<2$.
\endroster
\noindent
It may seem that the two cases may overlap. We shall see however that
they are actually distinct.

\medskip
{\smc Case (1).}
We fix $x_0>1$ with $g(x_0)\ge2$, and we choose a real $\alpha$ such that
$$
g(x_0)=x_0^\alpha+x_0^{-\alpha}.
\tag\DJ
$$
This is indeed possible since $g(x_0)\ge2$.

Starting from (\DJ), the two relations (\DG) and (\DH) 
allow us to compute $g(x)$ for all $x$ of the form
$x_{n,m}:=x_0^{2^{n}3^{-m}}$, $n,m\in\Bbb N_{\ge0}$. 
More precisely, assuming that we already know $g(x_{n,m})$
and that $g(x_{n,m})\ge2$ (which is certainly the case for $n=m=0$), 
by (\DG) and (\DH) we get 
$$\align 
g(x_{n+1,m})&=g^2(x_{n,m})-2,\\
g(x_{n,m})&=g^3(x_{n,m+1})-3g(x_{n,m+1}).
\endalign$$
In particular, we see that $g(x_{n+1,m})\ge2$. Moreover, 
since $g(x_{n,m})\ge2$, the second
equation determines the value $g(x_{n,m+1})$ {\it uniquely}, and it satisfies
$g(x_{n,m+1})\ge2$. 

A simple inductive argument --- the start of the induction being given by
(\DJ) --- then shows that, in
fact, the precise value of $g(x_{n,m})$ is given by (\DI).

We claim that the values $x_{n,m}=x_0^{2^{n}3^{-m}}$ are dense in the
positive real numbers. Clearly, it suffices to show that the
rational numbers $2^{n}3^{-m}$ are dense in the positive real numbers.
To see this, one takes the logarithm,
$$
\log_2\left(2^{n}3^{-m}\right)
=n-m\log_2(3).
$$
It is a well-known fact that
the sequence $(n\omega)_{n\ge0}$ is uniformly distributed modulo~1
if $\omega$ is irrational (cf\. \cite{\KuNiAA, Ex.~2.1}). 
Since $\log_2(3)$ is irrational, this implies our claim.

The function $g$ being continuous by assumption, this proves the
claim for Case~(1). 

\medskip
{\smc Case~(2).} We fix $x_0>1$ such that $g(x_0)<2$.
We claim that $\vert
g(x)\vert\le2$ for all $x>1$. Namely, it is certainly
not possible to find an $x_1>1$ with $g(x_1)\ge2$ because we would then
be in Case~(1), where we found that $g(x)\ge2$ for {\it all\/}
$x\in\Bbb R_{>0}$. Moreover, it is also not possible to 
find an $x_1>1$ with $g(x_1)\le-2$ because we would then obtain
$g(x_1^2)\ge2$ by (\DG), which would again put us in Case~(1). 

Now, using (\DG) iteratively, 
we will eventually find a positive integer $k$ such
that $g(x_0^{2^k})<0$. We recall that $g(1)=2$.
By the continuity of $g$ it then follows that 
there is a $y_0>1$ such that $g(y_0)=0$.
Without loss of generality, let $y_0$ be minimal with this property.
In particular, with this choice of $y_0$ we have $0< g(x)\le 2$
for all $x\in [1,y_0)$.

We now choose an $x_1\in(1,y_0)$. As we just said, we have
$0<g(x_1)\le 2$.
We may thus find an imaginary~$\alpha$ such that
$$
g(x_1)=x_1^\alpha+x_1^{-\alpha}=2\cos\big(i\alpha\log x_1\big).
\tag\DK
$$
Similarly to above, starting from this,
the two relations (\DG) and (\DH) allow us to
compute $g(x)$ for all $x\in[1,y_0]$ of the form $x_{n,m}=x_1^{2^{-n}3^m}$
with $n,m\in\Bbb N_{\ge0}$. 
More precisely, assuming that we already know $g(x_{n,m})$,
by (\DG) and (\DH) we get 
$$\align 
g(x_{n,m})&=g^2(x_{n+1,m})-2,\\
g(x_{n,m+1})&=g^3(x_{n,m})-3g(x_{n,m+1}),
\tag\DL
\endalign$$
as long as $x_{n+1,m}$ respectively $x_{n,m+1}$ is still in the
interval $[1,y_0)$. The important point here is that the first
equation determines the value $g(x_{n,m+1})$ {\it uniquely} as the
{\it positive} root of the equation since
we know that $g(x)>0$ for all $x\in[1,y_0)$.

Again, 
a simple inductive argument --- here starting from
(\DK) --- then shows that, in
fact, the precise value of $g(x_{n,m})$ is given by (\DI).
Finally, this may be extended to {\it all\/} $x_{n,m}$
(that is, without the restriction $x_{n,m}\in[1,y_0)$) by
applying (\DL) iteratively. Since also the numbers
$x_{n,m}=x_1^{2^{-n}3^m}$ are dense in the positive real numbers, this
establishes the claim for Case~(2), and thus completes the proof of the
lemma.\quad \quad \qed
\enddemo

\enddocument

\proclaim{Proposition}
Let $g:\Bbb R^{>0}\to\Bbb R$ be a continuous function satisfying the functional
equation
$$
g(x)g(y)=g(xy)+g(x/y).
\tag\DA
$$
Then either $g$ is identically zero or 
$$
g(x)=x^\alpha+x^{-\alpha},\quad \text{for }x\in\Bbb R_{>0},
$$
for a suitable $\alpha\in\Bbb R$.
\endproclaim

\demo{Proof}
Clearly, $g\equiv0$ solves (\DA). Hence, let us assume from now on
that $g\not\equiv0$.

By setting $y=1$ in (\DA), we see that $g(1)g(x)=2g(x)$. By our
assumption, this implies $g(1)=2$.

Next we set $x=1$ in (\DA). This leads to the relation
$$
g(y)=g(1/y),\quad \text{for }y\in\Bbb R_{>0}.
$$
This means that we may restrict our attention to the values $g(y)$
with $y\ge1$, which we shall do from now on.

We now claim that $g(y)\ge2$ for all $y\ge1$. Arguing by
contradiction, we assume that there is a $y_1>1$ with $g(y_1)<2$.
By setting $x=y$ in (\DA), we see that
$$
g^2(y)=g(y^2)+2,\quad \text{for }y\in\Bbb R_{>0}.
$$
or, equivalently,
$$
g(y^2)=g^2(y)-2,\quad \text{for }y\in\Bbb R_{>0}.
\tag\DB
$$
Either we have $g(y_1)<0$ or, by iteratively using (\DB), we may
find a $y_2>1$ with $g(y_2)<0$. Since $g(1)=2$ and $g$ is continuous,
the ?? theorem implies that there exists $y_3>1$ with $g(y_3)=0$.

We now take a $x_1$ strictly between $1$ and $y_3$ such that
$g(x_1)>0$, say $g(x_1)=a>0$. By substituting $y=y_3$ in (\DA),
we see that
$$
g(xy_3)=-g(x/y_3),
$$
or, equivalently,
$$
g(xy_3^2)=-g(x).
$$
Thus, we have $g(x_1y^{2n})=(-1)^na$ for $n=0,1,\dots$. 
Making use of the ?? theorem again, we obtain a sequence
$1<y_4<y_5<\cdots$ for which $g(y_i)=0$, $i=4,5,\dots$.

We now concentrate on $y_4$ and $y_5$ for which we have
$g(y_4)=g(y_5)=0$ and $g(x_2)\ne0$ for some $x_2$ between
$y_4$ and $y_5$.

Rewrite (\DA) in the form
$$
g\left(\sqrt{xy}\right)g\left(\sqrt{x/y}\right)=g(x)+g(y).
$$
If we choose $x$ and $y$ to be zeroes of the function $g$, then
at least one of $\sqrt{xy}$, $\sqrt{x/y}$, $\sqrt{y/x}$ is another
zero of~$g$. Using this argument iteratively, we will either the interval
$(1,y_4)$ or the interval $(y_4,y_5)$ (or both) with a dense set of
zeroes of~$g$, in contradiction with $g$ being continuous and
$g(1)=2$ respectively $g(x_2)\ne0$.

Thus, we indeed have $g(x)\ge2$ for all $x$.

Now fix $x_0>1$. We choose $\alpha\in\Bbb R$ so that
$$
g(x_0)=x_0^\alpha+x_0^{-\alpha}.
$$

Next, by setting $x=y^2$ in (\DA), we obtain
$$
g(y^2)g(y)=g(y^3)+g(y),\quad \text{for }y\in\Bbb R_{>0}.
$$
Recalling (\DB), we have
$$
g(y^2)=g^2(y)-2,\quad \text{for }y\in\Bbb R_{>0},
$$
and
$$
g(y^3)=g^3(y)-3g(y),\quad \text{for }y\in\Bbb R_{>0}.
$$
These two relations allow us to compute $g(x)$ for all $x$ of the form
$x_0^{2^{n_1}3^{-n_2}}$, $n_1,n_2\in\Bbb N_{\ge0}$. Since these numbers
are dense in $\Bbb R_{>0}$ (this boils down to the (stronger) fact that
$\log2/\log3$ is uniformly distributed modulo~1), we are
done.\quad \quad \qed
\enddemo

\proclaim{Proposition \TD}
Let $f$ and $g$ be continuous 
functions from $\Bbb R_{>0}$ to $\Bbb R$ satisfying the relations
$$\align 
f(a)g(x)&=f(ax)+f(a/x),
\tag\DA
\\
f(a)f(x)&=g(ax)-g(a/x)
\tag\DB
\endalign$$
for all positive real numbers $a$ and $x$. Then either $f$ is
identically zero and $g$ a constant function, or 
$$
f(x)=x^\alpha-x^{-\alpha}\quad \text{and}\quad 
g(x)=x^\alpha+x^{-\alpha}
$$
for all positive real numbers $x$, for a suitable real number $\alpha$.
\endproclaim

\demo{Proof}
Clearly, $f\equiv0$ and $g\equiv\text{const.}$ solves (\DA) and (\DB),
and nothing more general is possible for $g$ if $f\equiv0$.
Therefore, from now on, we may assume that $f$ is not identically
zero. 

We first find an $\alpha$ and a set $\{ x_{n,m} : n, m = 0, 1,
2, \dots  \}$  that is dense in $\Bbb R_{>0}$ such that
$$f(x)=x^\alpha-x^{-\alpha},$$
whenever $x=x_{n,m}$. By the continuity of $f$ this will imply that
this formula is true for any positive real number $x$. 

If we set $x=1$ in (\DB), then we obtain $f(a)f(1)=0$.
By our assumption on~$f$, this implies $f(1)=0$. Consequently,
if we set $a=1$ in (\DA), we obtain
$$
f(x)=-f(1/x), \quad \text{for }x\in\Bbb R_{>0}.
\tag\DC
$$
This relation implies in particular that it suffices to show
the claim of the proposition for $x\ge1$. Therefore, from now on we
concentrate on values $x$ in the interval $[1,\infty)$.

Now we multiply both sides of (\DB) by $f(y)$ and then use (\DA)
to obtain
$$
f(a)f(x)f(y)
=f(y) g( ax)
-f(y)g( a/x)
=f(axy)
+f(y/ax)
-f(ay/x)
-f(xy/a).
$$
Application of (\DC) with $x$ replaced by $y/ax$ turns this into
$$
f(a)f(x)f(y)
=f(axy)
-f(ax/y)
-f(ay/x)
-f(xy/a).
\tag\DD
$$
We choose $a=x=y$ in this relation to obtain
$$
f(x^3)
=f^3(x)+3f(x),\quad \text{for }x\in R_{>0}.
\tag\DE
$$
On the other hand, the choice of $a=x^3$ and $x=y$ in (\DD) leads to
$$
f(x^5)
=f^5(x)+5f^3(x)
+5f(x),\quad \text{for }x\in R_{>0},
\tag\DF
$$
if one makes again use of (\DC).

Let $x_0>1$. We may choose a real number $\alpha$ such that
$$
f(x_0)=x_0^\alpha-x_0^{-\alpha}.
$$
This is indeed possible since the right-hand side, as a function
in~$\alpha$, has the complete set of real numbers as range.

By (\DE), we then obtain
$$
f(x_0^3)=f^3(x_0)+3f(x_0)=x_0^{3\alpha}-x_0^{-3\alpha}
=\left(x_0^3\right)^{\alpha}-\left(x_0^3\right)^{-\alpha}.
$$
On the other hand, using (\DF) with $x=x_0^{1/5}$, we get the
polynomial equation 
$$
f(x_0)
=f^5\!\left(x_0^{1/5}\right)+5f^3\!\left(x_0^{1/5}\right)
+5f\!\left(x_0^{1/5}\right)$$
of degree~$5$ for $f\!\left(x_0^{1/5}\right)$. One can convince oneself 
that the polynomial $x^5+5x^3+5x$ is monotone increasing for real~$x$,
hence the above equation has a unique solution. This solution must be
$$
f\!\left(x_0^{1/5}\right)=x_0^{\alpha/5}-x_0^{-\alpha/5}
=\left(x_0^{1/5}\right)^{\alpha}-\left(x_0^{1/5}\right)^{-\alpha}.
$$

The above construction can be iterated. 
The upshot is that it shows that
$$
f\left(x_0^{3^{n}5^{-m}}\right)
=\left(x_0^{3^{n}5^{-m}}\right)^{\alpha}
-\left(x_0^{3^{n}5^{-m}}\right)^{-\alpha}
$$
for all non-negative integers $n$ and $m$.

We claim that the values $x_{n,m}=x_0^{3^{n}5^{-m}}$ are dense in the
positive real numbers. Clearly, it suffices to show that the
rational numbers $3^{n}5^{-m}$ are dense in the positive real numbers.
To see this, one takes the logarithm,
$$
\log_3\left(3^{n}5^{-m}\right)
=n-m\log_3(5).
$$
It is a well-known fact that
the sequence $(n\omega)_{n\ge0}$ is uniformly distributed modulo~1
if $\omega$ is irrational (cf\. \cite{\KuNiAA, Ex.~2.1}). 
Since $\log_3(5)$ is irrational, this implies our claim.

\medskip
The function $f$ being continuous, the above argument proves that 
$$
f(x)=x^\alpha-x^{-\alpha}
$$
for {\it all\/} positive real numbers $x\ge1$, and, by (\DC), for {\it
all\/} positive real numbers~$x$. One then substitutes this
in (\DA) to see that
$$
g(x)=x^\alpha+x^{-\alpha}
$$
for all positive real numbers~$x$.
This completes the proof of the proposition.\quad \quad \qed
\enddemo

\enddocument